\documentclass[12pt]{article}

\usepackage{amsmath,amssymb,amsthm,amscd,a4wide}

\begin{document}

\newcommand{\End}{{\rm{End}\ts}}
\newcommand{\Hom}{{\rm{Hom}}}
\newcommand{\Mat}{{\rm{Mat}}}
\newcommand{\ad}{{\rm{ad}\ts}}
\newcommand{\ch}{{\rm{ch}\ts}}
\newcommand{\chara}{{\rm{char}\ts}}
\newcommand{\diag}{ {\rm diag}}
\newcommand{\pr}{^{\tss\prime}}
\newcommand{\non}{\nonumber}
\newcommand{\wt}{\widetilde}
\newcommand{\wh}{\widehat}
\newcommand{\ot}{\otimes}
\newcommand{\la}{\lambda}
\newcommand{\ls}{\ts\lambda\ts}
\newcommand{\La}{\Lambda}
\newcommand{\De}{\Delta}
\newcommand{\al}{\alpha}
\newcommand{\be}{\beta}
\newcommand{\ga}{\gamma}
\newcommand{\Ga}{\Gamma}
\newcommand{\ep}{\epsilon}
\newcommand{\ka}{\kappa}
\newcommand{\vk}{\varkappa}
\newcommand{\vt}{\vartheta}
\newcommand{\si}{\sigma}
\newcommand{\vp}{\varphi}
\newcommand{\de}{\delta}
\newcommand{\ze}{\zeta}
\newcommand{\om}{\omega}
\newcommand{\ee}{\epsilon^{}}
\newcommand{\su}{s^{}}
\newcommand{\hra}{\hookrightarrow}
\newcommand{\ve}{\varepsilon}
\newcommand{\ts}{\,}
\newcommand{\vac}{\mathbf{1}}
\newcommand{\di}{\partial}
\newcommand{\qin}{q^{-1}}
\newcommand{\tss}{\hspace{1pt}}
\newcommand{\Sr}{ {\rm S}}
\newcommand{\U}{ {\rm U}}
\newcommand{\BL}{ {\overline L}}
\newcommand{\BE}{ {\overline E}}
\newcommand{\BP}{ {\overline P}}
\newcommand{\AAb}{\mathbb{A}\tss}
\newcommand{\CC}{\mathbb{C}\tss}
\newcommand{\KK}{\mathbb{K}\tss}
\newcommand{\QQ}{\mathbb{Q}\tss}
\newcommand{\SSb}{\mathbb{S}\tss}
\newcommand{\ZZ}{\mathbb{Z}\tss}
\newcommand{\X}{ {\rm X}}
\newcommand{\Y}{ {\rm Y}}
\newcommand{\Z}{{\rm Z}}
\newcommand{\Ac}{\mathcal{A}}
\newcommand{\Lc}{\mathcal{L}}
\newcommand{\Mc}{\mathcal{M}}
\newcommand{\Pc}{\mathcal{P}}
\newcommand{\Qc}{\mathcal{Q}}
\newcommand{\Tc}{\mathcal{T}}
\newcommand{\Sc}{\mathcal{S}}
\newcommand{\Bc}{\mathcal{B}}
\newcommand{\Ec}{\mathcal{E}}
\newcommand{\Fc}{\mathcal{F}}
\newcommand{\Hc}{\mathcal{H}}
\newcommand{\Uc}{\mathcal{U}}
\newcommand{\Vc}{\mathcal{V}}
\newcommand{\Wc}{\mathcal{W}}
\newcommand{\Yc}{\mathcal{Y}}
\newcommand{\Ar}{{\rm A}}
\newcommand{\Br}{{\rm B}}
\newcommand{\Ir}{{\rm I}}
\newcommand{\Fr}{{\rm F}}
\newcommand{\Jr}{{\rm J}}
\newcommand{\Or}{{\rm O}}
\newcommand{\GL}{{\rm GL}}
\newcommand{\Spr}{{\rm Sp}}
\newcommand{\Rr}{{\rm R}}
\newcommand{\Zr}{{\rm Z}}
\newcommand{\gl}{\mathfrak{gl}}
\newcommand{\middd}{{\rm mid}}
\newcommand{\ev}{{\rm ev}}
\newcommand{\Pf}{{\rm Pf}}
\newcommand{\Norm}{{\rm Norm\tss}}
\newcommand{\oa}{\mathfrak{o}}
\newcommand{\spa}{\mathfrak{sp}}
\newcommand{\osp}{\mathfrak{osp}}
\newcommand{\g}{\mathfrak{g}}
\newcommand{\h}{\mathfrak h}
\newcommand{\n}{\mathfrak n}
\newcommand{\z}{\mathfrak{z}}
\newcommand{\Zgot}{\mathfrak{Z}}
\newcommand{\p}{\mathfrak{p}}
\newcommand{\sll}{\mathfrak{sl}}
\newcommand{\agot}{\mathfrak{a}}
\newcommand{\qdet}{ {\rm qdet}\ts}
\newcommand{\Ber}{ {\rm Ber}\ts}
\newcommand{\HC}{ {\mathcal HC}}
\newcommand{\cdet}{ {\rm cdet}}
\newcommand{\tr}{ {\rm tr}}
\newcommand{\gr}{ {\rm gr}}
\newcommand{\str}{ {\rm str}}
\newcommand{\loc}{{\rm loc}}
\newcommand{\Gr}{{\rm G}}
\newcommand{\sgn}{ {\rm sgn}\ts}
\newcommand{\ba}{\bar{a}}
\newcommand{\bb}{\bar{b}}
\newcommand{\bi}{\bar{\imath}}
\newcommand{\bj}{\bar{\jmath}}
\newcommand{\bk}{\bar{k}}
\newcommand{\bl}{\bar{l}}
\newcommand{\hb}{\mathbf{h}}
\newcommand{\Sym}{\mathfrak S}
\newcommand{\fand}{\quad\text{and}\quad}
\newcommand{\Fand}{\qquad\text{and}\qquad}
\newcommand{\For}{\qquad\text{or}\qquad}
\newcommand{\OR}{\qquad\text{or}\qquad}

\renewcommand{\theequation}{\arabic{section}.\arabic{equation}}

\newtheorem{thm}{Theorem}[section]
\newtheorem{lem}[thm]{Lemma}
\newtheorem{prop}[thm]{Proposition}
\newtheorem{cor}[thm]{Corollary}
\newtheorem{conj}[thm]{Conjecture}
\newtheorem*{mthm}{Main Theorem}
\newtheorem*{mthma}{Theorem A}
\newtheorem*{mthmb}{Theorem B}

\theoremstyle{definition}
\newtheorem{defin}[thm]{Definition}

\theoremstyle{remark}
\newtheorem{remark}[thm]{Remark}
\newtheorem{example}[thm]{Example}

\newcommand{\bth}{\begin{thm}}
\renewcommand{\eth}{\end{thm}}
\newcommand{\bpr}{\begin{prop}}
\newcommand{\epr}{\end{prop}}
\newcommand{\ble}{\begin{lem}}
\newcommand{\ele}{\end{lem}}
\newcommand{\bco}{\begin{cor}}
\newcommand{\eco}{\end{cor}}
\newcommand{\bde}{\begin{defin}}
\newcommand{\ede}{\end{defin}}
\newcommand{\bex}{\begin{example}}
\newcommand{\eex}{\end{example}}
\newcommand{\bre}{\begin{remark}}
\newcommand{\ere}{\end{remark}}
\newcommand{\bcj}{\begin{conj}}
\newcommand{\ecj}{\end{conj}}

\newcommand{\bal}{\begin{aligned}}
\newcommand{\eal}{\end{aligned}}
\newcommand{\beq}{\begin{equation}}
\newcommand{\eeq}{\end{equation}}
\newcommand{\ben}{\begin{equation*}}
\newcommand{\een}{\end{equation*}}

\newcommand{\bpf}{\begin{proof}}
\newcommand{\epf}{\end{proof}}

\def\beql#1{\begin{equation}\label{#1}}

\title{\Large\bf Classical $\Wc$-algebras in types $A$, $B$, $C$, $D$ and $G$}

\author{{A. I. Molev\quad and\quad E. Ragoucy}}

\date{} % Start December 2013
\maketitle

\vspace{35 mm}

\begin{abstract}
We produce explicit generators of the classical $\Wc$-algebras associated
with the principal nilpotents in the simple Lie algebras of all classical types
and in the exceptional Lie algebra of type $G_2$. The generators are given
by determinant formulas in the context of the Poisson vertex algebras.
We also show that the images of the $\Wc$-algebra generators under
the Chevalley-type isomorphism
coincide with the elements defined via the corresponding Miura transformations.

% \medskip
%
% Mathematics Subject Classification 2010: 17B37, 17B67

\end{abstract}

%%%\vspace{5 mm}
%%%
%%%{\it Key words:}
%%%

\vspace{45 mm}

\noindent
School of Mathematics and Statistics\newline
University of Sydney,
NSW 2006, Australia\newline
alexander.molev@sydney.edu.au

\vspace{7 mm}

\noindent
Laboratoire de Physique Th\'{e}orique LAPTh,
CNRS and Universit\'{e} de Savoie\newline
BP 110, 74941 Annecy-le-Vieux Cedex, France\newline
eric.ragoucy@lapth.cnrs.fr

% \newpage
%
% \tableofcontents

\newpage

\section{Introduction}
\label{sec:int}
\setcounter{equation}{0}

Given a finite-dimensional simple Lie algebra $\g$ over $\CC$
and a nilpotent element $f\in\g$, the corresponding {\it classical
$\Wc$-algebra\/} $\Wc(\g,f)$ can be defined as a Poisson algebra of functions
on an infinite-dimensional manifold. The algebras $\Wc(\g,f)$ were
defined by Drinfeld and Sokolov~\cite{ds:la} in the case where $f$ is the
principal nilpotent and were used to introduce equations of the KdV type
for arbitrary simple Lie algebras.

Originally introduced in physics to investigate Toda systems \cite{ls:td} or
minimal representations of conformal field theories \cite{bbss:ev}, \cite{z:ia},
the $\Wc$-algebras were firstly studied in the context of constrained
Wess--Zumino--Witten models;
see \cite{bs:ws}, \cite{frrtw:hr},
and references therein for the physics literature.

In a recent work by De Sole, Kac and Valeri (see \cite{dskv:cw}  and review \cite{ds:cf})
the construction of Drinfeld and Sokolov (generalized to an arbitrary nilpotent
element $f$) was described in the framework
of Poisson vertex algebras. This description
was then applied to construct integrable hierarchies
of bi-Hamiltonian equations.

We follow the approach of \cite{dskv:cw} to produce explicit generators
of the classical
$\Wc$-algebra $\Wc(\g,f)$ for all simple Lie algebras $\g$ of types $A,B,C,D$ and $G$
for the case where $f$ is the principal nilpotent element of $\g$.
The generators are given as certain determinants of matrices formed by
elements of differential algebras associated with $\g$.
These formulas were suggested by the recent work \cite{am:eg}, where
the rectangular affine $\Wc$-algebras
were explicitly produced in type $A$.

We show that the images of the generators of the algebra $\Wc(\g,f)$
under a Chevalley-type projection coincide with the elements defined in \cite{ds:la}
via the Miura transformation. The projection yields an isomorphism
of Chevalley type to an algebra of polynomials defined as the intersection of the kernels
of screening operators. This provides a direct connection between two
presentations of the classical $\Wc$-algebras: one in the context of
the Poisson vertex algebras and the other as the Harish-Chandra image
of the Feigin--Frenkel center (i.e., the center of the affine vertex algebra
at the critical level); see \cite[Theorem~8.1.5]{f:lc}. Explicit generators of the center
were constructed in \cite{cm:ho}, \cite{ct:qs} and \cite{m:ff} for the Lie algebras $\g$
of all classical types. We make a connection between the Harish-Chandra images
of these elements ({\it loc. cit.} and \cite{mm:yc})
and the generators of $\Wc(\g,f)$ constructed via the determinant formulas.

The Chevalley-type isomorphism in the case of
the Lie algebra of type $G_2$ allows us to use the determinant
formula for the generators to produce explicit elements of $\Wc(\g,f)$
via a Miura transformation.

\medskip

This work was completed during the first author's visit to
the {\it Laboratoire
d'Annecy-le-Vieux de Physique Th\'{e}orique\/}. He
would like to thank the Laboratory for the warm hospitality.

\section{Definitions and preliminaries}
\label{sec:def}
\setcounter{equation}{0}

We will follow \cite{dskv:cw} to introduce the principal
classical $\Wc$-algebras in the context of Poisson vertex algebras.

\subsection{Differential algebras and $\la$-brackets}
\label{subsec:dif}

Let $\g$ be a vector space over $\CC$ and let
$X_1,\dots,X_d$ be a basis of $\g$.
Following \cite{dskv:cw}, consider the {\it differential algebra\/} $\Vc=\Vc(\g)$
which is defined as the {\it algebra
of differential polynomials\/}
in the variables $X_1,\dots,X_d$,
\ben
\Vc=\CC[X^{(r)}_1,\dots,X^{(r)}_d\ |\ r=0,1,2,\dots]\qquad 
\text{with}\quad X^{(0)}_i=X^{}_i,
\een
equipped with the derivation $\di$ defined by $\di\tss(X^{(r)}_i)=X^{(r+1)}_i$
for all $i=1,\dots,d$ and $r\geqslant 0$.

Suppose now that $\g$ is a simple (or reductive) Lie algebra and
fix a symmetric invariant bilinear form
$(.|.)$ on $\g$.
Introduce the $\la$-{\it bracket\/} on $\Vc$ as a linear map
\ben
\Vc\ot\Vc\to\CC[\la]\ot\Vc,\qquad a\ot b\mapsto \{a_{\ls}b\}.
\een
By definition, it is given by
\beql{defla}
\{X_{\ls}Y\}=[X,Y]+(X|Y)\tss\la\qquad\text{for}\quad X,Y\in\g,
\eeq
and extended to $\Vc$ by sesquilinearity $(a,b\in\Vc)$:
\ben
\{\di\tss a_{\ls}b\}=-\ls\tss\{a_{\ls}b\},\qquad
\{a_{\ls}\di\tss b\}=(\la+\di)\tss\{a_{\ls}b\},
\een
skewsymmetry
\ben
\{a_{\ls}b\}=-\{b_{\ts-\la-\di\ts}a\},
\een
and the Leibniz rule $(a,b,c\in\Vc)$:
\ben
\{a_{\ls}b\tss c\}=\{a_{\ls}b\}\tss c+\{a_{\ls}c\}\tss b.
\een
The $\la$-bracket defines the {\it affine Poisson vertex algebra\/} structure on
$\Vc$ \cite{dskv:cw}.

Now choose a Cartan subalgebra $\h$ of $\g$ and
a triangular decomposition $\g=\n_-\oplus\h\oplus \n_+$.
Set $\p=\n_-\oplus\h$ and define the projection map $\pi_{\p}:\g\to\p$
with the kernel $\n_+$.
From now on we will assume that $f\in\n_-$ is a principal nilpotent element of $\g$.
We regard $\Vc(\p)$ as a differential subalgebra of 
$\Vc$ and define the differential algebra homomorphism
\ben
\rho:\Vc\to \Vc(\p)
\een
by setting
\beql{rhopro}
\rho\tss(X)=\pi_{\p}(X)+(f\tss|X),\qquad X\in\g.
\eeq
The {\it classical $\Wc$-algebra $\Wc(\g,f)$} is defined by
\beql{defcwa}
\Wc(\g,f)=\{P\in\Vc(\p)\ |\ \rho\tss\{X_{\ls}P\}=0\quad\text{for all}\quad
X\in\n_+\}.
\eeq
By \cite[Lemma~3.2]{dskv:cw}, the classical $\Wc$-algebra $\Wc(\g,f)$
is a differential subalgebra of $\Vc(\p)$; moreover, it is a Poisson vertex algebra
equipped with the $\lambda$-bracket
\ben
\{a_{\ls}b\}_{\rho}=\rho\tss\{a_{\ls}b\},\qquad a,b\in \Wc(\g,f).
\een

We will need a sufficient condition for a family of elements of
the classical $\Wc$-algebra $\Wc(\g,f)$
to be its set of generators as a differential algebra (Proposition~\ref{prop:genwalg}).
Include the principal nilpotent element $f$ into an $\sll_2$-triple
$\{e,f,h\}$ with the standard commutation relations 
\ben
[e,f]=h,\qquad [h,e]=2\tss e,\qquad [h,f]=-2\tss f,
\een
and set $x=h/2$. The subspace $\p$ is assumed to be compatible
with the $\ad x$-eigenspace decomposition. In particular, $\p$ contains the
centralizer $\g^f$ of the element $f$ in $\g$.
Let $v_1,\dots,v_n$ be a basis
of $\g^f$ consisting of $\ad x$-eigenvectors.
If $b$ is an $\ad x$-eigenvector, we let $\de_x(b)$ denote its eigenvalue.
The following is a particular
case of \cite[Corollary~3.19]{dskv:cw}.

\bpr\label{prop:genwalg}
Suppose that $w_1,\dots,w_n$ is an arbitrary collection
of elements of $\Wc(\g,f)$ of the form $w_j=v_j+g_j$, where $g_j$
is the sum of products of $\ad x$-eigenvectors $b_i\in\p$,
\ben
g_j=\sum b_1^{(m_1)}\dots b_s^{(m_s)},\qquad
\sum_{i=1}^s\big(1-\de_x(b_i)+m_i\big)=1-\de_x(v_j),
\een
with $s+m_1+\dots+m_s>1$. Then $w_1,\dots,w_n$ are generators of the
differential algebra $\Wc(\g,f)$.
\qed
\epr

\subsection{Chevalley-type theorem}
\label{subsec:ctt}

Keeping the notation of the previous section, introduce the
standard Chevalley generators $e_i,h_i,f_i$ with $i=1,\dots,n$
of the simple Lie algebra $\g$ of rank $n$. The generators $h_i$ form
a basis of the Cartan subalgebra $\h$ of $\g$, while the $e_i$ and $f_i$
generate the respective nilpotent subalgebras $\n_+$ and $\n_-$.
Let $A=[a_{ij}]$
be the Cartan matrix of $\g$ so that the defining relations of $\g$ take the form
\begin{alignat}{2}
[e_i,f_j]&=\de_{ij}h_i,\qquad [h_i,h_j]&&=0,
\non\\
[h_i,e_j]&=a_{ij}\tss e_j,\qquad [h_i,f_j]&&=-a_{ij}\tss f_j,
\non
\end{alignat}
together with the Serre relations; see e.g. \cite{k:id}. There exists
a diagonal matrix $D=\diag\tss[\ep_1,\dots,\ep_n]$ with positive rational
entries such that the matrix $B=D^{-1}A$ is symmetric.
Following \cite[Ch.~2]{k:id}, normalize the symmetric invariant
bilinear form on $\g$ so that
\ben
(e_i|f_j)=\de_{ij}\tss \ep_i.
\een

We let $\De$ denote the root system of $\g$. For each root
$\al\in\De$ choose a nonzero root vector $e_{\al}$ so that
the set $\{e_{\al}\ |\ \al\in\De\}\cup\{h_i\ |\ i=1,\dots,n\}$ is a basis
of $\g$. The subset of positive roots will be denoted by $\De^+$. The
respective sets of elements $e_{\al}$ and  $e_{-\al}$ with $\al\in\De^+$
form bases of
$\n_+$ and $\n_-$. We also let $\al_1,\dots,\al_n$ denote the simple roots
so that
\ben
e_i=e_{\al_i},\qquad f_i=e_{-\al_i},\qquad i=1,\dots,n.
\een

The elements of the differential algebra $\Vc(\g)$
corresponding to the generators of $\g$
will be denoted
by
\ben
e^{(r)}_i=\di^{\tss r}(e_i),\qquad f^{(r)}_i=\di^{\tss r}(f_i),
\qquad h^{(r)}_i=\di^{\tss r}(h_i).
\een
Let
\beql{projh}
\phi:\Vc(\p)\to\Vc(\h)
\eeq
denote the homomorphism of differential algebras defined on the generators
as the projection $\p\to\h$ with the kernel $\n_-$. For each $i=1,\dots,n$
introduce the {\it screening operator\/}
\ben
V_i:\Vc(\h)\to \Vc(\h)
\een
defined by the formula
\beql{scr}
V_i=\sum_{r=0}^{\infty} V_{i\ts r}\sum_{j=1}^n a_{ji}\ts\frac{\di}{\di\tss h^{(r)}_j},
\eeq
where the coefficients $V_{i\ts r}$ are elements of $\Vc(\h)$ found by the
relation
\beql{genfv}
\sum_{r=0}^{\infty}\frac{V_{i\ts r}\ts z^r}{r!}=\exp\Big({-}\sum_{m=1}^{\infty}
\frac{h_i^{(m-1)}\tss z^m}{\ep_i\ts m!}\Big).
\eeq

The following statement is well-known; it
can be regarded as an analogue of the classical Chevalley theorem
providing a relationship between two presentations of the
classical $\Wc$-algebra; cf. \cite{dskv:cw},
\cite{ds:la} and \cite[Ch.~8]{f:lc}.

\bpr\label{prop:chev}
The restriction of the homomorphism $\phi$ to the classical $\Wc$-algebra
$\Wc(\g,f)$ yields an isomorphism
\beql{phiw}
\phi:\Wc(\g,f)\to \wt\Wc(\g,f),
\eeq
where $\wt\Wc(\g,f)$ is the subalgebra of $\Vc(\h)$ which consists of the elements
annihilated by all screening operators $V_i$,
\ben
\wt\Wc(\g,f)=\bigcap_{i=1}^n\ts \ker V_i.
\een
\epr

\bpf
We start by showing that for any $P\in \Wc(\g,f)$ its image $\phi(P)$
is annihilated by each operator $V_i$. It will be convenient to work
with an equivalent affine version of
the differential algebra $\Vc(\g)$; see \cite[Sec.~2.7]{k:va}.
Recall that the affine Kac--Moody algebra $\wh\g$
is defined as the central
extension
\beql{km}
\wh\g=\g\tss[t,t^{-1}]\oplus\CC K,
\eeq
where $\g[t,t^{-1}]$ is the Lie algebra of Laurent
polynomials in $t$ with coefficients in $\g$. For any $r\in\ZZ$ and $X\in\g$
we set $X[r]=X\ts t^r$. The commutation relations of the Lie algebra $\wh\g$
have the form
\ben
\big[X[r],Y[s]\big]=[X,Y][r+s]+r\ts\de_{r,-s}(X|Y)\ts K,
\qquad X, Y\in\g,
\een
and the element $K$ is central in $\wh\g$. Consider the quotient $\Sr(\wh\g)/\Ir$
of the symmetric algebra $\Sr(\wh\g)$ by its ideal $\Ir$
generated by the subspace $\g[t]$ and the element
$K-1$. This quotient can be identified
with the symmetric algebra $\Sr\big(t^{-1}\g[t^{-1}]\big)$,
as a vector space.
We will identify the differential algebras $\Vc(\g)\cong \Sr\big(t^{-1}\g[t^{-1}]\big)$
via the isomorphism
\beql{isaff}
X^{(r)}\mapsto r!\ts X[-r-1],\qquad X\in\g,\quad r\geqslant 0,
\eeq
so that the derivation $\di$ will correspond to the derivation $-d/d\tss t$.
Similarly, we will identify the differential
algebras $\Vc(\p)\cong \Sr\big(t^{-1}\p[t^{-1}]\big)$.

By definition~\eqref{defcwa}, if an element $P\in\Vc(\p)$ belongs to the
subalgebra $\Wc(\g,f)$, then $\rho\tss\{{e_i}^{}_{\ls}P\}=0$ for all $i=1,\dots,n$.
Now observe that, regarding $P$ as an element of the $\g[t]$-module
$\Sr(\wh\g)/\Ir\cong \Sr\big(t^{-1}\g[t^{-1}]\big)$,
we can write
\ben
\{{e_i}^{}_{\ls}P\}=\sum_{r=0}^{\infty}\frac{\la^r}{r!}\ts e_i[r]\ts P.
\een
We have the following relations in $\wh\g$,
\ben
\bal
\big[e_i[r], f_i[-s-1]\big]&=h_i[r-s-1]+r\ts\de_{r,\tss s+1}\ts\ep_i\tss K ,\\[0.3em]
\big[e_i[r], h_j[-s-1]\big]&=-a_{ji}\ts e_i[r-s-1].
\eal
\een
Moreover, for each positive root $\al\ne\al_i$ we also have
\ben
\big[e_i[r], e_{-\al}[-s-1]\big]=c_i(\al)\ts e_{-\al+\al_i}[r-s-1]
\een
for a certain constant $c_i(\al)$, if $\al-\al_i$ is a root;
otherwise the commutator is zero.
Hence, recalling the definition \eqref{rhopro} of the homomorphism $\rho$,
we can conclude that the condition that $P$ belongs to the
subalgebra $\Wc(\g,f)$ implies the relations
\beql{eip}
\wh e_i[r]\ts P=0\qquad\text{for all}\quad i=1,\dots,n\fand r\geqslant 0,
\eeq
where $\wh e_i[r]$ is the operator on $\Sr\big(t^{-1}\p[t^{-1}]\big)$
given by
\ben
\bal
\wh e_i[r]&=\sum_{s=r}^{\infty}h_i[r-s-1]\ts\frac{\di}{\di f_i[-s-1]}
+\ep_i\ts r\ts\frac{\di}{\di f_i[-r]}
-\sum_{j=1}^n a_{ji} \ts\frac{\di}{\di h_j[-r-1]}\\
{}&+\sum_{\al\in\De^+,\ts \al\ne\al_i}\ts
\sum_{s=r}^{\infty}c_i(\al)\ts e_{-\al+\al_i}[r-s-1]\ts \frac{\di}{\di\ts e_{-\al}[-s-1]},
\eal
\een
and $e_{-\al+\al_i}$ is understood as being equal to zero, if $\al-\al_i$ is not a root.
Denote the generating function in $z$ introduced in \eqref{genfv} by $V_i(z)$
and replace $h_i^{(m-1)}/(m-1)!$ with $h_i[-m]$ for $m\geqslant 1$
in accordance with \eqref{isaff}.
We have the relation for its derivative,
\ben
V'_i(z)=V_i(z)\Big({-}\sum_{m=1}^{\infty}
\frac{h_i[-m]\tss z^{m-1}}{\ep_i}\Big).
\een
Taking the coefficient of $z^{p-1}$ with $p\geqslant 1$ we get the relations
\beql{cance}
\ep_i\ts p \ts \frac{V_{i\ts p}}{p!}+\sum_{r=0}^{p-1}\ts \frac{V_{i\ts r}}{r!}
\ts h_i[r-p]=0.
\eeq
By \eqref{eip} the element $P$ has the property
\beql{modic}
\sum_{r=0}^{\infty}\frac{V_{i\ts r}}{r!}\ts \wh e_i[r]\ts P=0.
\eeq
Note that by \eqref{cance}
all the differentiations $\di/\di f_i[-s-1]$
with $s\geqslant 0$  will cancel in the expansion of the left hand side
of \eqref{modic}. Moreover, the elements of the form $e_{-\al+\al_i}[r-s-1]$
occurring in the expansion of $\wh e_i[r]$ will vanish under the projection
\eqref{projh}. Therefore, \eqref{modic} implies that the image
$\phi(P)$ with respect to this projection satisfies the relation
\ben
\sum_{r=0}^{\infty}\frac{V_{i\ts r}}{r!}\ts
\sum_{j=1}^n a_{ji} \ts\frac{\di}{\di h_j[-r-1]}\ts\phi(P)=0
\een
which is equivalent to
\ben
\sum_{r=0}^{\infty} V_{i\ts r}\sum_{j=1}^n a_{ji}\ts\frac{\di}{\di\tss h^{(r)}_j}
\ts \phi(P)=0,
\een
that is, $V_i\ts \phi(P)=0$, as claimed. This shows
that the restriction of the homomorphism $\phi$ to the classical $\Wc$-algebra
$\Wc(\g,f)$ provides the homomorphism \eqref{phiw}.

It is clear from the definition that $\phi$ is a differential algebra homomorphism.
It is known that $\Wc(\g,f)$ is an algebra of differential polynomials
in $n$ variables (a proof of this claim for an arbitrary nilpotent element $f$
is given in \cite[Theorem~~3.14]{dskv:cw}).
Therefore, the proof of the proposition
can be completed by verifying that the images of the generators
of $\Wc(\g,f)$ under the homomorphism $\phi$ are differential algebra
generators of $\wt\Wc(\g,f)$. We will omit this verification for the general case.
In types $A$, $B$, $C$, $D$ and $G$ this will follow easily from
the explicit construction of generators of $\Wc(\g,f)$ given
in the next sections.
\epf

\section{Generators of $\Wc(\gl_n,f)$}
\label{sec:gen}
\setcounter{equation}{0}

We will consider square matrices of the form
\beql{matramod}
A=\begin{bmatrix}a_{11}&a_{12}&0&0&\dots&0\\
                 a_{21}&a_{22}&a_{23}&0&\dots&0\\
                 a_{31}&a_{32}&a_{33}&a_{34}&\dots&0\\
                 \dots&\dots&\dots&\dots&\dots &\dots \\
                             a_{n-1\tss 1}&a_{n-1\tss 2}&a_{n-1\tss 3}
                             &\dots&\dots&a_{n-1\tss n}\\
                             a_{n1}&a_{n2}&a_{n3}&\dots&\dots&a_{nn}
                \end{bmatrix}
\eeq
such that the entries $a_{ij}$ belong
to a certain ring.
It is well-known
and can be easily verified that
even if the ring
is noncommutative, the column-determinant and row-determinant of the matrix $A$
coincide, provided that the entries $a_{12}$, $a_{23},\dots, a_{n-1\tss n}$
belong to the center of the ring. We will assume that this condition holds,
and define the (noncommutative) determinant of $A$ by setting
\beql{det}
\det A=\sum_{\si\in\Sym_n}\sgn\si\cdot a_{\si(1)\ts 1}\dots a_{\si(n)\ts n}
=\sum_{\si\in\Sym_n}\sgn\si\cdot a_{1\ts\si(1)}\dots a_{n\ts \si(n)}.
\eeq

A particular case which we will often use below
is where the ring contains the identity
and the matrix $A$ has the form
\beql{matra}
A=\begin{bmatrix}a_{11}&1&0&0&\dots&0\\
                 a_{21}&a_{22}&1&0&\dots&0\\
                 a_{31}&a_{32}&a_{33}&1&\dots&0\\
                 \dots&\dots&\dots&\dots&\dots &\dots \\
                             a_{n-1\tss 1}&a_{n-1\tss 2}&a_{n-1\tss 3}&\dots&\dots&1\\
                             a_{n1}&a_{n2}&a_{n3}&\dots&\dots&a_{nn}
                \end{bmatrix}.
\eeq
One easily derives the following explicit formula
for the determinant of this matrix:
\beql{detaone}
\det A=(-1)^{n-1}\ts a^{}_{n\tss 1}+\sum_{1\leqslant i_1<\dots<i_k<n}
(-1)^{n-k-1}\ts a^{}_{i_1 1}\ts a^{}_{i_2\ts i_1+1}\ts a^{}_{i_3\ts i_2+1}\dots
a^{}_{n\ts i_k+1},
\eeq
where $k$ runs over the set of values $1,\dots,n-1$.

Returning to more general matrices \eqref{matramod},
denote by $D_i$ (resp., $\overline D_i$) the determinant of the $i\times i$ submatrix
of $A$ corresponding to the first (resp., last) $i$ rows and columns. We suppose that
$D_0=\overline D_0=1$.

\ble\label{lem:deex}
Fix $p\in\{0,1,\dots,n\}$. Then for the determinant of the matrix \eqref{matramod}
we have
\ben
\det A=D_p\ts \overline D_{n-p}+
\sum_{j=1}^{p}\sum_{i=p+1}^n\ts(-1)^{j+i}\ts D_{j-1}\ts \wt a^{}_{i\tss j}\ts \overline D_{n-i},
\een
where
\beql{modme}
\wt a^{}_{i\tss j}=a^{}_{i\tss j}\ts
a_{j\ts j+1}\ts a_{j+1\ts j+2}\dots a_{i-1\ts i}\qquad\text{for}\quad i> j.
\eeq
\ele

\bpf
The formula follows easily from the definition of the determinant \eqref{det}.
In the applications which we consider below, the central elements
$a_{12}$, $a_{23},\dots, a_{n-1\tss n}$ turn out to be invertible.
Then the lemma is reduced to the particular case of matrices \eqref{matra}
and it is immediate from the explicit formula \eqref{detaone}.
\epf

\subsection{Determinant-type generators}

Take $\g$ to be the general linear Lie algebra $\gl_n$ with its standard basis $E_{ij}$,
$i,j=1,\dots,n$. The elements $E_{11},\dots,E_{nn}$ span a Cartan subalgebra
of $\gl_n$ which we will denote by $\h$.
The respective subsets of basis elements $E_{ij}$ with $i<j$ and $i>j$ span
the nilpotent subalgebras $\n_+$ and $\n_-$. The subalgebra $\p=\n_-\oplus\h$ is then spanned
by the elements $E_{ij}$ with $i\geqslant j$. Take the principal nilpotent element $f$
in the form
\ben
f=E_{2\tss 1}+E_{3\tss 2}+\dots+E_{n\tss n-1}.
\een
For the $\sll_2$-triple $\{e,f,h\}$ take
\ben
e=\sum_{i=1}^{n-1} i(n-i)\ts E_{i\ts i+1}\Fand
h=\sum_{i=1}^{n} (n-2\tss i+1)\ts E_{i\tss i}.
\een

We will be working with the algebra of differential operators $\Vc(\p)\ot\CC[\di]$,
where the commutation relations are given by
\ben
\di\ts E^{(r)}_{ij}- E^{(r)}_{ij}\ts \di=E^{\tss(r+1)}_{ij}.
\een
In other words, $\di$ will be regarded as a generator of this algebra rather
than the derivation on $\Vc(\p)$. For any element $g\in\Vc(\p)$
and any nonnegative integer $r$ the element $g^{(r)}=\di^{\tss r}(g)$
coincides with
the constant term of the differential operator $\di^{\tss r} g$ so that
\beql{grdi}
g^{(r)}=\di^{\tss r} g\ts 1,
\eeq
assuming that $\di\ts 1=0$.

The invariant symmetric bilinear form on $\gl_n$
is defined by
\beql{invforma}
(X|Y)=\tr\ts XY,\qquad X,Y\in\gl_n,
\eeq
where $X$ and $Y$ are understood as $n\times n$ matrices over $\CC$.

Consider the determinant \eqref{det} of the matrix with entries in $\Vc(\p)\ot\CC[\di]$,
\beql{deta}
\det
\begin{bmatrix}\di+E_{11}&1&0&0&\dots&0\\
                 E_{21}&\di+E_{22}&1&0&\dots&0\\
                 E_{31}&E_{32}&\di+E_{33}&1&\dots&0\\
                 \dots&\dots&\dots&\dots&\dots& \dots \\
                             E_{n-1\tss 1}&E_{n-1\tss 2}&E_{n-1\tss 3}&\dots&\dots&1\\
                             E_{n1}&E_{n2}&E_{n3}&\dots&\dots&\di+E_{nn}
                \end{bmatrix}
                =\di^{\ts n}+w_1\ts \di^{\ts n-1}+\dots+w_n.
\eeq

\bth\label{thm:gln}
All elements $w_1,\dots,w_n$ belong to the classical $\Wc$-algebra
$\Wc(\gl_n,f)$. Moreover, the elements $w^{(r)}_1,\dots,w^{(r)}_n$
with $r=0,1,\dots$ are algebraically independent
and generate the algebra $\Wc(\gl_n,f)$.
\eth

\bpf
Denote the determinant on the left hand side of \eqref{deta} by $D_n$.
For each $1\leqslant k<n$ we will identify the Lie algebra $\gl_k$ with the natural
subalgebra of $\gl_n$ so that the determinant $D_k$ will also be regarded
as an element of the algebra $\Vc(\p)\ot\CC[\di]$. We also set $D_0=1$.
We will be proving
the first part of the theorem by induction on $n$. The induction base is trivial
and for any $1\leqslant i<n-1$ the induction hypothesis implies
\beql{rhoei}
\rho\tss\{E_{i\ts i+1}\ts^{}_{\ls} D_k\}=0
\eeq
for all $k=i+1,\dots,n$. Moreover, relation \eqref{rhoei} clearly holds
for the values $k=1,\dots,i-1$ as well, while
\ben
\rho\tss\{E_{i\ts i+1}\ts^{}_{\la} D_i\}
=-\rho\tss\big(D^+_{i-1}\ts E_{i\ts i+1}\big)
=-D^+_{i-1},
\een
where for each polynomial $P=P(\di)\in \Vc(\p)\ot\CC[\di]$
we denote by $P^+$
the polynomial $P(\di+\la)$ regarded
as an element of $\CC[\la]\ot\Vc(\p)\ot\CC[\di]$.
By Lemma~\ref{lem:deex} (with $p=n-1$) we have the relation
\begin{multline}
D_n=D_{n-1}\ts (\di+E_{nn})-D_{n-2}\ts E_{n\tss n-1}+D_{n-3}\ts E_{n\tss n-2}\\[0.3em]
{}+\dots+(-1)^{n-2} \ts D_1\tss E_{n\ts 2}+(-1)^{n-1} \ts D_0\tss E_{n\tss 1}.
\label{dnexpa}
\end{multline}
Hence, using the properties of the $\la$-bracket we get
\ben
\rho\tss\{E_{i\ts i+1}\ts^{}_{\la} D_n\}=
(-1)^{n-i+1}\ts
D^+_{i-1}\ts E_{n\tss i+1}+(-1)^{n-i}\ts
D^+_{i-1}\ts E_{n\tss i+1}=0
\een
for all $i=1,\dots,n-2$. Furthermore,
\begin{multline}
\rho\tss\{E_{n-1\ts n}\ts^{}_{\la} D_n\}=D^+_{n-1}
-D^+_{n-2}\ts  (\di+E_{nn})-D^+_{n-2}\ts  (E_{n-1\ts n-1}-E_{n\tss n}+\la)\\[0.3em]
+D^+_{n-3}\ts E_{n-1\ts n-2}
{}+\dots+(-1)^{n-2} \ts D^+_1\tss E_{n-1\ts 2}+(-1)^{n-1} \ts D^+_0\tss E_{n-1\tss 1}
\non
\end{multline}
which is zero due to relation \eqref{dnexpa} applied to the determinant $D^+_{n-1}$
instead of $D_n$. Since $\rho\tss\{E_{i\ts i+1}\ts^{}_{\la} D_n\}=0$
for all $i=1,\dots,n-1$, we may conclude that
$\rho\tss\{X\tss^{}_{\la} D_n\}=0$ for all $X\in\n_+$ so that
all elements $w_1,\dots,w_n$ belong to the subalgebra $\Wc(\gl_n,f)$
of $\Vc(\p)$.

To prove the second part of the theorem we apply Proposition~\ref{prop:genwalg}.
Note that the powers of the $n\times n$ matrix $f$ form a basis of the centralizer
$\gl_n^{\tss f}$
so that for $j=1,\dots,n$ we can take $v_j$ to be equal, up to a sign, to $f^{j-1}$.
Then the condition of Proposition~\ref{prop:genwalg} holds for
the elements $w_1,\dots,w_n$ thus implying that they are generators
of the differential algebra $\Wc(\gl_n,f)$.

Finally, the images of the elements $w_k$ under the homomorphism \eqref{phiw}
are the elements $\wt w_k\in\Vc(\h)$ found from the relation
\ben
(\di+E_{11})\dots (\di+E_{n\tss n})=\di^{\ts n}+\wt w_1\ts \di^{\ts n-1}+\dots+\wt w_n.
\een
In the notation of Sec.~\ref{subsec:ctt} we have
\ben
h^{(r)}_j=E^{(r)}_{j\tss j}-E^{(r)}_{j+1\ts j+1},\qquad j=1,\dots,n-1.
\een
The Cartan matrix is of the size $(n-1)\times(n-1)$,
\ben
A=\begin{bmatrix}2&-1\phantom{-}&0&\dots&0&0\\
                 -1\phantom{-}&2&-1\phantom{-}&\dots&0&0\\
                 0&-1\phantom{-}&2&\dots&0&0\\
                 \dots&\dots&\dots&\ddots&\dots &\dots \\
                             0&0&0&\dots&2&-1\phantom{-}\\
                             0&0&0&\dots&-1\phantom{-}&2
                \end{bmatrix}
\een
so that $a_{ii}=2$ for $i=1,\dots,n-1$
and $a_{i\ts i+1}=a_{i+1\ts i}=-1$ for $i=1,\dots,n-2$, while all other
entries are zero. Hence, regarding $\Vc(\h)$ as the algebra of polynomials
in the variables $E^{(r)}_{i\tss i}$ with $i=1,\dots,n$ and $r=0,1,\dots$, we get
\ben
\sum_{j=1}^{n-1} a_{ji}\ts\frac{\di}{\di\tss h^{(r)}_j}=
\frac{\di}{\di\tss E^{(r)}_{i\tss i}}-\frac{\di}{\di\tss E^{(r)}_{i+1\ts i+1}}.
\een
Therefore, the screening operators \eqref{scr} take the form
\ben
V_i=\sum_{r=0}^{\infty} V_{i\ts r}
\Bigg(\frac{\di}{\di\tss E^{(r)}_{i\tss i}}
-\frac{\di}{\di\tss E^{(r)}_{i+1\ts i+1}}\Bigg),\qquad i=1,\dots,n-1,
\een
where the coefficients $V_{i\ts r}$ are found by the
relation
\ben
\sum_{r=0}^{\infty}\frac{V_{i\ts r}\ts z^r}{r!}=\exp\Big({-}\sum_{m=1}^{\infty}
\frac{E^{(m-1)}_{i\tss i}-E^{(m-1)}_{i+1\ts i+1}}{m!}\ts z^m\Big).
\een
The differential algebra $\wt\Wc(\gl_n,f)$ consists of the polynomials
in the variables $E^{(r)}_{i\tss i}$, which are annihilated by all
operators $V_i$. It is easy to verify directly
that the elements $\wt w_1,\dots,\wt w_n$ belong to $\wt\Wc(\gl_n,f)$; cf. \cite{mm:yc}.
Moreover, the elements $\wt w^{\tss(r)}_1,\dots,\wt w^{\tss(r)}_n$
with $r$ running over nonnegative integers
are algebraically independent
generators of the algebra $\wt\Wc(\gl_n,f)$; see \cite[Ch.~8]{f:lc}.
Hence, the generators $w^{\tss(r)}_1,\dots,w^{\tss(r)}_n$ of $\Wc(\gl_n,f)$
are also algebraically independent.
\epf

The injective homomorphism $\Wc(\gl_n,f)\hra \Vc(\h)$
taking $w_i$ to $\wt w_i$ constructed
in the proof of Theorem~\ref{thm:gln}
can be shown to respect the $\la$-brackets. This fact
is known as the {\it Kupershmidt--Wilson theorem\/}, and the embedding
is called the {\it Miura transformation\/};
see e.g. \cite{a:tf}, \cite[Ch.~8]{f:lc},
\cite{gd:fh} and \cite{kw:ml}.
Note also that the arguments used in the proof imply that
the homomorphism \eqref{phiw}
is bijective in the case $\g=\gl_n$; see
the proof of Proposition~\ref{prop:chev}.

This embedding can be used
to calculate $\la$-brackets in the classical $\Wc$-algebra inside $\Vc(\h)$,
as we illustrate below; cf. \cite{v:cw}.
The $\la$-bracket on $\Vc(\h)$
is determined by the relations
\ben
\{{E_{ii}}_{\ts\la\ts} E_{jj}\}=\de_{i\tss j}\ts\la.
\een
Let us set $x_i=\di+E_{ii}$ for $i=1,\dots,n$ and let $u$ be a variable
commuting with the $x_i$. Setting $C=E_{11}+\dots+E_{nn}$,
we find
\ben
\{C_{\ts\la\ts} (u+x_1)\dots(u+x_n)\}=(u+x^+_1)\dots(u+x^+_n)-(u+x_1)\dots(u+x_n),
\een
where $x^+_i=\di+\la+E_{ii}$. Similarly,
consider the quadratic element
\ben
P=-\frac12\ts \sum_{i=1}^n E^2_{i\tss i}
-\sum_{i=1}^n (n-i)\ts E_{i\tss i}.
\een
Then
\ben
\{P_{\ts\la\ts} (u+x_1)\dots(u+x_n)\}=u\ts (u+x^+_1)\dots(u+x^+_n)
-(\di+n\tss\la+u)(u+x_1)\dots(u+x_n),
\een
which follows by a straightforward calculation based on
the properties of the $\la$-brackets.
Defining the coefficients $\wt e_m$ of the polynomial in $u$,
by
\ben
(u+x_1)\dots(u+x_n)=\sum_{m=0}^n \wt e_m\ts u^{n-m}
\een
we get
\ben
\{P_{\ts\la\ts} \wt e_m\}=\wt e^{\ts +}_{m+1}-\wt e_{m+1}-(\di+n\tss\la)\tss \wt e_m,
\een
where $\wt e^{\ts +}_{m+1}$ is obtained from $\wt e_{m+1}$ by  replacing $\di$ with
$\di+\la$. Note that the $\wt e_m$
coincide with the images of the respective differential operators $e_m$
under the isomorphism \eqref{phiw}; see \eqref{emexp} below.

\subsection{MacMahon theorem}
\label{subsec:mmt}

We now return to an arbitrary matrix $A$ of the form \eqref{matra}
with entries in
a certain ring with the identity. Let $t$ be a variable commuting
with elements of the ring.
Define elements $e_m$ of the ring by the expansion
\beql{dete}
\det(1+tA)=\sum_{m=0}^{n} e_m\tss t^m
\eeq
and denote this polynomial in $t$ by $e(t)$.

\ble\label{lem:em}
The elements $e_m$ are found by
\ben
e_m=\sum_{s=1}^m\ts\sum_{\ i_k\geqslant j_k\ts\ts\text{and}\ts\ts\ts i_k<j_{k+1}}
(-1)^{m-s}\ts a_{i_1\tss j_1}\dots a_{i_s\tss j_s},
\een
where the second sum is taken
over the indices $i_1,\dots,i_s$ and
$j_1,\dots,j_s$ such that $i_k\geqslant j_{k}$ for $k=1,\dots,s$ and
$i_k<j_{k+1}$ for $k=1,\dots,s-1$ satisfying the condition
\beql{sumjik}
\sum_{k=1}^s (j_k-i_k)=m-s.
\eeq
\ele

\bpf
For any $1\leqslant k\leqslant n$ we will use the notation $e^{\{k\}}_m$
to indicate the elements $e_m$ associated with the submatrix of $A$
corresponding to the first $k$ rows and columns.
The obvious analogue
of \eqref{dnexpa} for the determinant $\det(1+tA)$ gives
the recurrence relation
\begin{multline}
e^{\{n\}}_m=e^{\{n-1\}}_m+e^{\{n-1\}}_{m-1}\ts a_{n\tss n}
-e^{\{n-2\}}_{m-2}\ts a_{n\ts n-1}\\[0.3em]
{}+\dots+(-1)^{m-2}\ts
e^{\{n-m+1\}}_{1}\ts a_{n\ts n-m+2}+(-1)^{m-1}\ts
a_{n\ts n-m+1}
\non
\end{multline}
which implies the desired formula.
\epf

For an arbitrary matrix $A$ of the form \eqref{matra} over a ring, we also define
the family of elements $h_m$ of the ring by setting $h_0=1$ and
\ben
h_m=\sum_{s=1}^{\infty}\ts\sum_{\ i_k\geqslant j_k\ts\ts\text{and}\ts\ts\ts i_k
\geqslant j_{k+1}}
a_{i_1\tss j_1}\dots a_{i_s\tss j_s},\qquad m\geqslant 1,
\een
where the second sum is taken
over the indices $i_1,\dots,i_s$ and
$j_1,\dots,j_s$ such that $i_k\geqslant j_{k}$ for $k=1,\dots,s$ and
$i_k\geqslant j_{k+1}$ for $k=1,\dots,s-1$
satisfying \eqref{sumjik}.
Combine
these elements into the formal series
\ben
h(t)=\sum_{m=0}^{\infty} h_m\tss t^m.
\een

The following identity is a version of the MacMahon Master Theorem
for noncommutative matrices of the form \eqref{matra}.

\bpr\label{prop:mmta}
We have the identity
\ben
h(t)\ts e(-t)=1.
\een
\epr

\bpf
We need to verify that for any $m\geqslant 1$ we have the relation
\beql{reche}
h_m-h_{m-1}\ts e_1+h_{m-2}\ts e_2+\dots+(-1)^n h_{m-n}\ts e_n=0,
\eeq
where we assume that $h_k=0$ for $k<0$. The expansion of the left hand side
as a linear combination of monomials in the entries of the matrix $A$
will contain monomials of the form
\ben
a_{i_1\tss j_1}\dots a_{i_p\tss j_p}\ts a_{i_{p+1}\tss j_{p+1}}\dots a_{i_{p+r}\tss j_{p+r}}
\een
such that $i_k\geqslant j_{k+1}$ for $k=1,\dots,p$
and $i_k< j_{k+1}$ for $k=p+1,\dots,p+r-1$. Such a monomial will occur
twice
in the expansion with the opposite signs.
Indeed, it will occur
in the expansion of $h_{a}\ts e_b$ for two pairs of indices $(a,b)$.
In the first pair,
\ben
a=\sum_{k=1}^p (j_k-i_k)+p,\qquad b=\sum_{k=p+1}^{p+r} (j_k-i_k)+r
\een
so that the sign of the monomial in the expansion is $(-1)^r$. In the second pair,
\ben
a=\sum_{k=1}^{p+1} (j_k-i_k)+p+1,\qquad b=\sum_{k=p+2}^{p+r} (j_k-i_k)+r-1,
\een
and the sign of the monomial in the expansion is $(-1)^{r-1}$.
\epf

\bre\label{rem:nonsf}
One can regard $e(t)$ and $h(t)$ as specializations of the generating
series of the noncommutative elementary and complete symmetric functions,
respectively; see \cite[Sec.~3]{gkllrt:ns}. Explicit formulas for
the specializations of the power sums symmetric functions and ribbon Schur
functions are also easy to write down in terms of the entries of $A$.
\qed
\ere

\subsection{Permanent-type generators}
\label{subsec:perm}

Now we specialize the matrix $A$ given in \eqref{matra} and take
\ben
a_{i\tss j}=\de_{ij}\tss\di+E_{i\tss j},\qquad i\geqslant j,
\een
so that the entries belong to the algebra $\Vc(\p)\ot\CC[\di]$.
Consider the elements $e_m$ and $h_m$ associated with this specialization
of the matrix $A$ and write them as differential operators,
\begin{align}
e_m&=e_{m\ts 0}+e_{m\ts 1}\ts\di+\dots+e_{m\ts m}\ts\di^{\ts m},
\label{emexp}\\[0.2em]
h_m&=h_{m\ts 0}+h_{m\ts 1}\ts\di+\dots+h_{m\ts m}\ts\di^{\ts m}.
\label{hmexp}
\end{align}
In particular, the constant terms are found by the formulas
\ben
e_{m\ts 0}=\sum_{s=1}^m\ts\sum_{\ i_k\geqslant j_k\ts\ts\text{and}\ts\ts\ts i_k<j_{k+1}}
(-1)^{m-s}\ts\big(\de_{i_1\tss j_1}\tss
\di+E_{i_1\tss j_1}\big)\dots \big(\de_{i_s\tss j_s}\tss
\di+E_{i_s\tss j_s}\big)\ts 1,
\een
where the second sum is taken
over the indices $i_1,\dots,i_s$ and
$j_1,\dots,j_s$ such that $i_k\geqslant j_{k}$ for $k=1,\dots,s$ and
$i_k<j_{k+1}$ for $k=1,\dots,s-1$
satisfying \eqref{sumjik}; and
\ben
h_{m\ts 0}=\sum_{s=1}^{\infty}\ts\sum_{\ i_k\geqslant j_k\ts\ts
\text{and}\ts\ts\ts i_k\geqslant j_{k+1}}\tss
\big(\de_{i_1\tss j_1}\tss\di+E_{i_1\tss j_1}\big)\dots \big(\de_{i_s\tss j_s}\tss
\di+E_{i_s\tss j_s}\big)\ts 1,
\een
where the second sum is taken
over the indices $i_1,\dots,i_s$ and
$j_1,\dots,j_s$ such that $i_k\geqslant j_{k}$ for $k=1,\dots,s$ and
$i_k\geqslant j_{k+1}$ for $k=1,\dots,s-1$
satisfying \eqref{sumjik}.

\bpr\label{prop:ewcoin}
We have
\beql{emiw}
e_{m\ts i}=\binom{n-m+i}{i}\ts w_{m-i}\qquad\text{for all}\quad
0\leqslant i\leqslant m\leqslant n,
\eeq
where the elements $w_1,\dots,w_n$ are introduced in \eqref{deta}.
\epr

\bpf
Consider the determinant $D_n=D_n(\di)$ defined in \eqref{deta}.
We have
\ben
\det(1+tA)=t^{\tss n} D_n(\di+t^{-1}).
\een
Hence
\ben
\sum_{m=0}^n\sum_{i=0}^m\ts e_{m\ts i}\ts \di^{\tss i}\ts t^m=
\sum_{k=0}^n\ts w_{k}\ts t^{\tss k}\ts (1+t\tss\di)^{n-k}
\een
which implies \eqref{emiw}.
\epf

In particular, $e_{m\ts 0}=w_{m}$ and so the family
$e^{(r)}_{1\tss 0},\dots,e^{(r)}_{n\tss 0}$
with $r=0,1,\dots$
is algebraically independent
and generates the algebra $\Wc(\gl_n,f)$.

\bco\label{cor:genh}
All elements $h_{m\ts i}$ belong to $\Wc(\gl_n,f)$.
Moreover, the family
$h^{(r)}_{1\tss 0},\dots,h^{(r)}_{n\tss 0}$
with $r=0,1,\dots$
is algebraically independent
and generates the algebra $\Wc(\gl_n,f)$.
\eco

\bpf
By Proposition~\ref{prop:ewcoin},
each $e_m$ is a linear combination
of differential operators whose coefficients
belong to $\Wc(\gl_n,f)$. Furthermore, Proposition~\ref{prop:mmta} implies that
each $h_m$ is a differential operator with coefficients
in $\Wc(\gl_n,f)$. This proves the first part of the corollary.

Furthermore, it is easy to verify (see also \cite[Sec.~4.1]{gkllrt:ns}) that
the recurrence relation \eqref{reche}
is equivalent to
\ben
h_m=\sum_{i_1+\dots+i_k=m}
(-1)^{m-k}\ts e^{}_{i_1}\dots e^{}_{i_k},
\een
summed over $k$-tuples of positive integers $(i_1,\dots,i_k)$ with
$k=1,\dots,m$. Together with \eqref{emiw}
this implies an explicit expression for $h_{m\tss 0}$ in terms of
the generators $e^{(r)}_{1\tss 0},\dots,e^{(r)}_{n\tss 0}$,
\ben
h_{m\tss 0}=(-1)^{m-1}\ts e_{m\tss 0}+\sum\ts \text{const}\cdot
\ts e^{(r_1)}_{j_1\ts 0}\dots e^{(r_k)}_{j_k\ts 0},
\een
summed over $k$-tuples of positive integers $(j_1,\dots,j_k)$
and $k$-tuples of nonnegative integers $(r_1,\dots,r_k)$ such that
$j_1+\dots+j_k+r_1+\dots+r_k=m$, with certain coefficients, where for
the summands with $k=1$ we have $r_1\geqslant 1$. This shows that
the elements $h^{(r)}_{1\tss 0},\dots,h^{(r)}_{n\tss 0}$
with $r=0,1,\dots$
are algebraically independent
and generate the algebra $\Wc(\gl_n,f)$.
\epf

\section{Generators of $\Wc(\oa_{2n+1},f)$}
\label{sec:genb}
\setcounter{equation}{0}

Given a positive integer $N$, we will use the involution on the set
$\{1,\dots,N\}$ defined by $i\mapsto i^{\tss\prime}=N-i+1$.
The Lie subalgebra
of $\gl_N$ spanned by the elements
\beql{fijo}
F_{i\tss j}=E_{i\tss j}-E_{j^{\tss\prime}\tss i^{\tss\prime}},\qquad i,j=1,\dots,N,
\eeq
is
the {\it orthogonal Lie algebra\/} $\oa_N$.
The simple Lie algebras of type $B_n$ correspond to the odd values
$N=2n+1$, while the simple Lie algebras of type $D_n$
correspond to the even values $N=2n$.
In both cases, we have
the commutation relations
\beql{on}
[F_{i\tss j}, F_{k\tss l}]=\de_{k\tss j}\tss F_{i\tss l}-\de_{i\tss l}\tss F_{k\tss l}
-\de_{k\tss i^{\tss\prime}}\tss F_{j^{\tss\prime}\tss l}
+\de_{j^{\tss\prime}\tss l}\tss F_{k\tss i^{\tss\prime}},
\qquad i,j,k,l\in\{1,\dots,N\}.
\eeq
Note also the symmetry relation
\beql{symrel}
F_{i\tss j}+F_{j^{\tss\prime}\tss i^{\tss\prime}}=0, \qquad i,j\in\{1,\dots,N\}.
\eeq

The elements $F_{11},\dots,F_{nn}$ span a Cartan subalgebra
of $\oa_N$ which we will denote by $\h$.
The respective subsets of elements $F_{ij}$ with $i<j$ and $i>j$ span
the nilpotent subalgebras $\n_+$ and $\n_-$. The subalgebra $\p=\n_-\oplus\h$ is then spanned
by the elements $F_{ij}$ with $i\geqslant j$.

As with the case of $\gl_n$,
we will be working with the algebra of differential operators $\Vc(\p)\ot\CC[\di]$,
where the commutation relations are given by
\beql{dfpol}
\di\ts F^{(r)}_{ij}- F^{(r)}_{ij}\ts \di=F^{\tss(r+1)}_{ij}.
\eeq
For any element $g\in\Vc(\p)$
and any nonnegative integer $r$ the element $g^{(r)}$
coincides with the constant term of the differential operator $\di^{\tss r} g$
as in \eqref{grdi}.

In this section we will work with the odd orthogonal Lie algebras $\oa_{2n+1}$.
Take the principal nilpotent element $f\in\oa_{2n+1}$
in the form
\beql{fbn}
f=F_{2\tss 1}+F_{3\tss 2}+\dots+F_{n+1\tss n}.
\eeq
The $\sll_2$-triple is now formed by the elements $\{e,f,h\}$ with
\ben
e=\sum_{i=1}^{n} i\tss(2\tss n-i+1)\ts F_{i\ts i+1}\Fand
h=2\ts\sum_{i=1}^{n} (n-i+1)\ts F_{i\tss i}.
\een
The invariant symmetric bilinear form on $\oa_{2n+1}$
is defined by
\beql{invformb}
(X|Y)=\frac12\ts \tr\ts XY,\qquad X,Y\in\oa_{2n+1},
\eeq
where $X$ and $Y$ are understood as matrices over $\CC$
which are skew-symmetric with respect to the antidiagonal.

Consider the determinant \eqref{det} of the matrix with entries in $\Vc(\p)\ot\CC[\di]$,
\ben
\det
\begin{bmatrix}\di+F_{1\tss 1}&1&\dots&0&0&0&\dots&0&0\\
                 F_{2\tss 1}&\di+F_{2\tss 2}&\dots&0&0&0&\dots&0&0\\
                 \dots&\dots&\ddots&\dots&\dots&\dots&\dots&\dots&\dots\\
                 F_{n\tss 1}&F_{n\tss 2}&\dots&\di+F_{n\tss n}&1&0&\dots&0&0\\
                 F_{n+1\tss 1}&F_{n+1\tss 2}&\dots&F_{n+1\ts n}&\di&-1\phantom{-}&\dots&0&0 \\
                 F_{n^{\tss\prime}\tss 1}&F_{n^{\tss\prime}\tss 2}
                            &\dots&0&F_{n^{\tss\prime}\tss n+1}
                            &\di+F_{n^{\tss\prime}\tss n^{\tss\prime}}&\dots&0&0\\
                 \dots&\dots&\dots&\dots&\dots&\dots&\ddots&\dots&\dots\\
                 F_{2^{\tss\prime}\tss 1}&0&\dots&\dots&F_{2^{\tss\prime}\tss n+1}
                          &\dots&\dots&\di+F_{2^{\tss\prime}\tss 2^{\tss\prime}}&-1\phantom{-}\\
                 0&F_{1'\tss 2}&\dots&\dots&F_{1'\tss n+1}
                           &\dots&\dots&F_{1'\tss 2^{\tss\prime}}&\di+F_{1'\tss 1'}
                \end{bmatrix}
\een
so that the $(i,j)$ entry of the matrix is $\de_{ij}\tss\di+F_{i\tss j}$ for $i\geqslant j$,
the $(i,i+1)$ entry is $1$ for $i\leqslant n$ and $-1$ for $i>n$, while the remaining entries
are zero. The determinant has the form
\beql{diffo}
\di^{\ts 2n+1}+w_2\ts \di^{\ts 2n-1}+w_3\ts \di^{\ts 2n-2}+\dots+w_{2n+1},
\qquad w_i\in\Vc(\p).
\eeq

\bth\label{thm:b}
All elements $w_2,\dots,w_{2n+1}$ belong to the classical $\Wc$-algebra
$\Wc(\oa_{2n+1},f)$. Moreover, the elements $w^{(r)}_2,w^{(r)}_4,\dots,w^{(r)}_{2n}$
with $r=0,1,\dots$ are algebraically independent
and generate the algebra $\Wc(\oa_{2n+1},f)$.
\eth

\bpf
The argument is similar to the proof of Theorem~\ref{thm:gln}.
Denote the determinant by $D$ and let $D_i$ (resp., $\overline D_i$) denote
the $i\times i$ minor
corresponding to the first (resp., last) $i$ rows and columns. We suppose that
$D_0=\overline D_0=1$. Lemma~\ref{lem:deex} implies the expansion
\ben
D=D_n\ts \di \ts\overline D_n+\sum_{j,k=1}^{n+1}(-1)^{n-j+1}\ts
D_{j-1}\ts F_{k'\tss j}\ts \overline D_{k-1}.
\een

To prove the first part of the theorem, note that the elements $F_{ij}$
with $1\leqslant i,j\leqslant n$ span a subalgebra of $\oa_{2n+1}$ isomorphic
to $\gl_n$. Hence, due to Theorem~\ref{thm:gln}, if $1\leqslant i\leqslant n-1$ then
\beql{rhoeib}
\rho\tss\{F_{i\ts i+1}\ts^{}_{\la}\ts D_k\}=0\Fand
\rho\tss\{F_{i\ts i+1}\ts^{}_{\la}\ts \overline D_k\}=0
\eeq
for all $1\leqslant k\leqslant n$ with $k\ne i$.
Moreover,
\ben
\rho\tss\{F_{i\ts i+1}\ts^{}_{\la}\ts D_i\}
=-D^+_{i-1}\Fand \rho\tss\{F_{i\ts i+1}\ts^{}_{\la}\ts \overline D_i\}
=\overline D_{i-1},
\een
where, as before, $P^+=P(\di+\la)$ for any polynomial
$P=P(\di)\in \Vc(\p)\ot\CC[\di]$. Hence, for any $k\in\{1,\dots,n+1\}$
and $k\ne i,i+1$
we have
\ben
\rho\tss\{F_{i\ts i+1}\ts^{}_{\la}\ts D_{i}\ts F_{k'\tss i+1}\ts \overline D_{k-1}
-D_{i-1}\ts F_{k'\tss i}\ts \overline D_{k-1}\}=0.
\een
Similarly, for any $j\in\{1,\dots,n+1\}$
and $j\ne i,i+1$
we have
\ben
\rho\tss\{F_{i\ts i+1}\ts^{}_{\la}\ts D_{j-1}\ts F_{i'\tss j}\ts \overline D_{i-1}
+D_{j-1}\ts F_{(i+1)'\tss j}\ts \overline D_{i}\}=0
\een
and
\ben
\rho\tss\{F_{i\ts i+1}\ts^{}_{\la}\ts D_{i-1}\ts F_{(i+1)'\tss i}\ts \overline D_{i}
-D_{i}\ts F_{i'\tss i+1}\ts \overline D_{i-1}\}=0.
\een
These relations imply $\rho\tss\{F_{i\ts i+1}\ts^{}_{\la}\ts D\}=0$. Finally,
performing similar calculations we get
\begin{multline}
\rho\tss\{F_{n\ts n+1}\ts^{}_{\la}\ts D\}=D^+_{n}\ts (\di+\la) \ts\overline D_{n-1}
+\sum_{k=1}^{n-1}D^+_n\ts F_{k'\tss n'}\ts \overline D_{k-1}-
D^+_n\ts (F_{n\tss n}+\la)\ts \overline D_{n-1}
\\
-D^+_{n-1}\ts \di \ts\overline D_n
+\sum_{j=1}^{n-1}(-1)^{n-j+1}\ts D^+_{j-1}\ts F_{n\ts j}\ts \overline D_{n}-
D^+_{n-1}\ts (F_{n\tss n}+\la)\ts \overline D_{n}.
\non
\end{multline}
Applying Lemma~\ref{lem:deex} to the determinants $D^+_{n}$ and $\overline D_{n}$,
we get the relations
\beql{dplusex}
D^+_{n}=D^+_{n-1}\ts (\di+\la+F_{n\tss n})
+\sum_{j=1}^{n-1}(-1)^{n-j}\ts D^+_{j-1}\ts F_{n\ts j}
\eeq
and
\beql{dbarplusex}
\overline D_{n}=(\di+F_{n'\tss n'})\ts \overline D_{n-1}
+\sum_{k=1}^{n-1}\ts F_{k'\tss n'}\ts \overline D_{k-1}
\eeq
which imply that $\rho\tss\{F_{n\ts n+1}\ts^{}_{\la}\ts D\}=0$.
This shows that
all elements $w_2,\dots,w_{2n+1}$ belong to the subalgebra $\Wc(\oa_{2n+1},f)$
of $\Vc(\p)$.

Now use Proposition~\ref{prop:genwalg}.
The odd powers of the matrix $f$ form a basis of the centralizer
$\oa_{2n+1}^{\tss f}$
so that for $j=1,\dots,n$ we can take $v_j$ to be equal, up to a sign, to $f^{2\tss j-1}$.
Then the condition of Proposition~\ref{prop:genwalg} will hold for the family
of elements $w_2,w_4,\dots,w_{2n}$,
thus implying that they are generators
of the differential algebra $\Wc(\oa_{2n+1},f)$.

The images of the elements $w_k$ under the homomorphism \eqref{phiw}
are the elements $\wt w_k\in\Vc(\h)$ found from the relation
\begin{multline}
(\di+F_{11})\dots (\di+F_{n\tss n})\ts\di\ts
(\di+F_{n'\tss n'}) \dots (\di+F_{1'\tss 1'})\\[0.3em]
=\di^{\ts 2n+1}+\wt w_2\ts \di^{\ts 2n-1}+\wt w_3\ts \di^{\ts 2n-2}+\dots+\wt w_{2n+1}.
\non
\end{multline}
In the notation of Sec.~\ref{subsec:ctt} we have
\ben
h^{(r)}_j=F^{(r)}_{j\tss j}-F^{(r)}_{j+1\ts j+1},\quad j=1,\dots,n-1,\Fand
h^{(r)}_n=2\tss F^{(r)}_{n\tss n}.
\een
The Cartan matrix is of the size $n\times n$,
\ben
A=\begin{bmatrix}2&-1\phantom{-}&0&\dots&0&0\\
                 -1\phantom{-}&2&-1\phantom{-}&\dots&0&0\\
                 0&-1\phantom{-}&2&\dots&0&0\\
                 \dots&\dots&\dots&\ddots&\dots &\dots \\
                             0&0&0&\dots&2&-1\phantom{-}\\
                             0&0&0&\dots&-2\phantom{-}&2
                \end{bmatrix}
\een
so that $a_{ii}=2$ for $i=1,\dots,n$
and $a_{i\ts i+1}=a_{i+1\ts i}=a_{n-1\ts n}=-1$ for $i=1,\dots,n-2$, while
$a_{n\ts n-1}=-2$ and
all other
entries are zero. The entries of the diagonal matrix $D=\diag\tss[\ep_1,\dots,\ep_n]$
are found by
\ben
\ep_1=\dots=\ep_{n-1}=1\Fand \ep_n=2.
\een
Hence, regarding $\Vc(\h)$ as the algebra of polynomials
in the variables $F^{(r)}_{i\tss i}$ with $i=1,\dots,n$ and $r=0,1,\dots$, we get
\ben
\sum_{j=1}^{n} a_{ji}\ts\frac{\di}{\di\tss h^{(r)}_j}=
\frac{\di}{\di\tss F^{(r)}_{i\tss i}}-\frac{\di}{\di\tss F^{(r)}_{i+1\ts i+1}},
\qquad i=1,\dots,n-1,
\een
and
\ben
\sum_{j=1}^{n} a_{j\tss n}\ts\frac{\di}{\di\tss h^{(r)}_j}=
\frac{\di}{\di\tss F^{(r)}_{n\tss n}}.
\een
Therefore, the screening operators \eqref{scr} take the form
\ben
V_i=\sum_{r=0}^{\infty} V_{i\ts r}
\Bigg(\frac{\di}{\di\tss F^{(r)}_{i\tss i}}
-\frac{\di}{\di\tss F^{(r)}_{i+1\ts i+1}}\Bigg),\qquad i=1,\dots,n-1,
\een
and
\ben
V_n=\sum_{r=0}^{\infty} V_{n\ts r}\ts\frac{\di}{\di\tss F^{(r)}_{n\tss n}},
\een
where the coefficients $V_{i\ts r}$ are found by the
relations
\ben
\sum_{r=0}^{\infty}\frac{V_{i\ts r}\ts z^r}{r!}=\exp\Big({-}\sum_{m=1}^{\infty}
\frac{F^{(m-1)}_{i\tss i}-F^{(m-1)}_{i+1\ts i+1}}{m!}\ts z^m\Big),\qquad i=1,\dots,n-1,
\een
and
\ben
\sum_{r=0}^{\infty}\frac{V_{n\ts r}\ts z^r}{r!}=\exp\Big({-}\sum_{m=1}^{\infty}
\frac{F^{(m-1)}_{n\tss n}}{m!}\ts z^m\Big).
\een
The differential algebra $\wt\Wc(\oa_{2n+1},f)$ consists of the polynomials
in the variables $F^{(r)}_{i\tss i}$, which are annihilated by all
operators $V_i$. It is easy to verify directly
that the elements $\wt w_k$ belong to $\wt\Wc(\oa_{2n+1},f)$;
cf. \cite[Sec.~4.2]{mm:yc}.
Furthermore, the elements $\wt w^{\tss(r)}_2,\wt w^{\tss(r)}_4,\dots,
\wt w^{\tss(r)}_{2n}$
with $r$ running over nonnegative integers
are algebraically independent
generators of the algebra $\wt\Wc(\oa_{2n+1},f)$; see \cite[Ch.~8]{f:lc}.
Hence, the generators $w^{\tss(r)}_2,w^{\tss(r)}_4,\dots,
w^{\tss(r)}_{2n}$ of $\Wc(\oa_{2n+1},f)$
are also algebraically independent.
\epf

The injective homomorphism $\Wc(\oa_{2n+1},f)\hra \Vc(\h)$
taking $w_i$ to $\wt w_i$
is known as the {\it Miura transformation\/} in type $B$ \cite{ds:la};
see also \cite[Ch.~8]{f:lc}. Note that by the arguments used in
the proof of Theorem~\ref{thm:b},
the homomorphism \eqref{phiw}
is bijective in the case $\g=\oa_{2n+1}$.

\bre\label{rem:fold}
There is an alternative way to prove the first part of Theorem~\ref{thm:b}
based on the idea of {\it folding\/} which was already used in the work \cite{frs:fw};
see also \cite[Ch.~2]{v:cw}.
To outline the argument, let $\theta$ denote the involutive automorphism of
the differential algebra $\Vc=\Vc(\gl_{2n+1})$
defined by
\ben
\theta: E^{(r)}_{ij}\mapsto -E^{(r)}_{j\pr i\pr}.
\een
We have the direct sum decomposition
\ben
\gl_{2n+1}=\oa_{2n+1}\oplus \gl^{\tss-}_{2n+1},
\een
where $\oa_{2n+1}$ is identified with the fixed point subalgebra under $\theta$, while
\ben
\gl^{\tss-}_{2n+1}=\{X\in\gl_{2n+1}\ |\ \theta(X)=-X\}
\een
is the subspace of anti-invariants.
With the identification $\Vc=\Sr\big(t^{-1}\gl_{2n+1}[t^{-1}]\big)$
used in Sec.~\ref{subsec:ctt},
for the subalgebra of $\theta$-invariants in $\Vc$
we have
\ben
\Vc^{\ts\theta}=\Sr\Big(t^{-1}\oa_{2n+1}[t^{-1}]\oplus
\Sr^2\big(t^{-1}\gl^{\tss-}_{2n+1}[t^{-1}]\big)\Big).
\een
This subalgebra is stable under the $\la$-bracket so that
its restriction defines
a $\la$-bracket on $\Vc^{\ts\theta}$.
Let $\Jr$ be the ideal of $\Vc^{\ts\theta}$ generated by the
subspace $\Sr^2\big(t^{-1}\gl^{\tss-}_{2n+1}[t^{-1}]\big)$.
For any element $P\in \Vc^{\ts\theta}$
we have $\{P^{}_{\ts\la\ts}\Jr\}\subset \Jr$.
The quotient space $\Vc^{\ts\theta}/\Jr$ is naturally identified
with the differential algebra $\Vc(\oa_{2n+1})$.
This quotient is equipped with a $\la$-bracket induced from
that of $\Vc^{\ts\theta}$. The resulting bracket
on $\Vc(\oa_{2n+1})$ is then obtained as the
{\it folding\/} of the $\la$-bracket on $\Vc(\gl_{2n+1})$.
It coincides with the $\la$-bracket
defined in \eqref{defla}.

For any element $P\in\Vc$ denote by $\text{\rm pr\tss}(P)$ its projection
to the subspace $\Vc^{\ts\theta}$ along the second summand
in the decomposition
\ben
\Vc=\Vc^{\ts\theta}\oplus t^{-1}\gl^{\tss-}_{2n+1}[t^{-1}]\ts\Vc^{\ts\theta}.
\een
Furthermore, let $\wt P\in\Vc(\oa_{2n+1})$ be the image of $\text{\rm pr\tss}(P)$
under the natural
epimorphism
\ben
\Vc^{\ts\theta}\to \Vc^{\ts\theta}/\Jr\cong \Vc(\oa_{2n+1}).
\een
Set
\ben
\wt f=E_{2\tss 1}+E_{3\tss 2}+\dots+E_{n+1\tss n}-E_{n+2\tss n+1}-\dots-E_{2n+1\tss 2n}
\in\gl_{2n+1}.
\een
Theorem~\ref{thm:gln} remains valid with the modified determinant
(of the size $(2n+1)\times(2n+1)$)
obtained from \eqref{deta} by changing the sign of the entries $1$
in the last $n$ columns.
Note that $\wt f$ coincides with the element $f$ defined in \eqref{fbn}.
The corresponding coefficients of the differential
operator will belong to the classical $\Wc$-algebra $\Wc(\gl_{2n+1},\wt f)$.
Moreover, the images of these coefficients under the map $P\mapsto\wt P$
belong to the classical $\Wc$-algebra $\Wc(\oa_{2n+1},f)$.
Upon an appropriate re-normalizing of the form \eqref{invforma}
to match \eqref{invformb}, we recover the
corresponding coefficients of the differential operator
\eqref{diffo} thus proving that they belong to $\Wc(\oa_{2n+1},f)$.

The same argument can be applied to the case of the
classical $\Wc$-algebra $\Wc(\spa_{2n},f)$; cf. Sec.~\ref{sec:genc} below.
Some modifications are required for the case of $\Wc(\g_{2},f)$
as we need to deal with a third order automorphism of $\oa_8$;
see Sec.~\ref{sec:geng} for the details. Note, however, that such a folding
argument is not applicable for the reduction in type $D$, i.e., for the
reduction from $\gl_{2n}$ to $\oa_{2n}$. The reason, which was already
pointed out in \cite{ds:la}, is the fact that the principal nilpotent element
$f\in\oa_{2n}$ cannot be lifted to a principal nilpotent
of the Lie algebra $\gl_{2n}$ containing $\oa_{2n}$.
\qed
\ere

Now we use Proposition~\ref{prop:mmta} to produce another family of generators
of $\Wc(\oa_{2n+1},f)$; cf. Corollary~\ref{cor:genh}.
Observe that the determinant yielding the differential operator \eqref{diffo}
has the form \eqref{matramod}. As in \eqref{modme},
for any indices
$2n+1\geqslant i\geqslant j\geqslant 1$ set
\beql{wtab}
\wt a_{i\tss j}=\begin{cases}\de_{ij}\tss\di+F_{i\tss j}\qquad&\text{if}
\quad n+1\geqslant i\geqslant j,\\[0.2em]
(-1)^{i-n-1}\ts\big(\de_{ij}\tss\di+F_{i\tss j}\big)\qquad&\text{if}
\quad i\geqslant n+1\geqslant j,\\[0.2em]
(-1)^{i-j}\ts\big(\de_{ij}\tss\di+F_{i\tss j}\big)\qquad&\text{if}
\quad i\geqslant j\geqslant n+1.
\end{cases}
\eeq
By analogy with Sec.~\ref{subsec:perm}
introduce elements $e_{m\ts 0},h_{m\ts 0}\in \Vc(\p)$ as constant terms
of the differential operators,
\beql{emob}
e_{m\ts 0}=\sum_{s=1}^m\ts\sum_{\ i_k\geqslant j_k\ts\ts\text{and}\ts\ts\ts i_k<j_{k+1}}
(-1)^{m-s}\ts\ts
\wt a_{i_1\tss j_1}\dots \wt a_{i_s\tss j_s}\ts 1,
\eeq
where the second sum is taken
over the indices $i_1,\dots,i_s$ and
$j_1,\dots,j_s$ such that $i_k\geqslant j_{k}$ for $k=1,\dots,s$ and
$i_k<j_{k+1}$ for $k=1,\dots,s-1$
satisfying \eqref{sumjik}; and
\beql{hmob}
h_{m\ts 0}=\sum_{s=1}^{\infty}\ts\sum_{\ i_k\geqslant j_k\ts\ts
\text{and}\ts\ts\ts i_k\geqslant j_{k+1}}\tss
\wt a_{i_1\tss j_1}\dots \wt a_{i_s\tss j_s}\ts 1,
\eeq
where the second sum is taken
over the indices $i_1,\dots,i_s$ and
$j_1,\dots,j_s$ such that $i_k\geqslant j_{k}$ for $k=1,\dots,s$ and
$i_k\geqslant j_{k+1}$ for $k=1,\dots,s-1$
satisfying \eqref{sumjik}. Proposition~\ref{prop:ewcoin} implies
$e_{m\ts 0}=w_{m}$ for all $m=2,\dots,2n+1$.
Hence, by Theorem~\ref{thm:b} the family
$e^{(r)}_{2\ts 0},e^{(r)}_{4\ts 0},\dots,e^{(r)}_{2n\ts 0}$
with $r=0,1,\dots$
is algebraically independent
and generates the algebra $\Wc(\oa_{2n+1},f)$.

\bco\label{cor:genhb}
All elements $h_{m\ts 0}$ belong to the algebra $\Wc(\oa_{2n+1},f)$.
Moreover, the family
$h^{(r)}_{2\ts 0},h^{(r)}_{4\ts 0},\dots,h^{(r)}_{2n\ts 0}$
with $r=0,1,\dots$
is algebraically independent
and generates $\Wc(\oa_{2n+1},f)$.
\eco

\bpf
The proof is essentially the same as that of Corollary~\ref{cor:genh};
the only difference is that the elements $e_{2k+1\ts 0}$ (resp., $h_{2k+1\ts 0}$)
with $k=0,1,\dots$ are expressible in terms of the elements $e_{2k\ts 0}$
(resp., $h_{2k\ts 0}$).
\epf

\section{Generators of $\Wc(\spa_{2n},f)$}
\label{sec:genc}
\setcounter{equation}{0}

As in Sec.~\ref{sec:genb}, for any
positive integer $n$, we use the involution on the set
$\{1,\dots,2n\}$ defined by $i\mapsto i^{\tss\prime}=2n-i+1$.
The Lie subalgebra
of $\gl_{2n}$ spanned by the elements
\beql{fijsp}
F_{ij}=E_{ij}-\ve_i\tss\ve_j\tss E_{j'i'},\qquad i,j=1,\dots,2n,
\eeq
is
the {\it symplectic Lie algebra\/} $\spa_{2n}$, where
we set $\ve_i=1$ for $i=1,\dots,n$ and
$\ve_i=-1$ for $i=n+1,\dots,2n$.
This is a
simple Lie algebra of type $C_n$.
We have
the commutation relations
\beql{spn}
[F_{ij}, F_{kl}]=\de_{kj}\tss F_{il}-\de_{il}\tss F_{kl}
-\ve_i\tss\ve_j\tss\big(\de_{ki'}\tss F_{j'l}-\de_{j'l}\tss F_{ki'}\big)
\eeq
for all $i,j,k,l\in\{1,\dots,2n\}$.
Note also the symmetry relation
\ben
F_{ij}+\ve_i\tss\ve_j\tss F_{j'i'}=0, \qquad i,j\in\{1,\dots,2n\}.
\een
The elements $F_{11},\dots,F_{nn}$ span a Cartan subalgebra
of $\spa_{2n}$ which we will denote by $\h$.
The respective subsets of elements $F_{ij}$ with $i<j$ and $i>j$ span
the nilpotent subalgebras $\n_+$ and $\n_-$. The subalgebra $\p=\n_-\oplus\h$ is then spanned
by the elements $F_{ij}$ with $i\geqslant j$.

As before,
we will work with the algebra of differential operators $\Vc(\p)\ot\CC[\di]$,
where the commutation relations are given by
\ben
\di\ts F^{(r)}_{ij}- F^{(r)}_{ij}\ts \di=F^{\tss(r+1)}_{ij}.
\een
For any element $g\in\Vc(\p)$
and any nonnegative integer $r$ the element $g^{(r)}$
coincides with
the constant term of the differential operator $\di^{\tss r} g$
as in \eqref{grdi}.

Take the principal nilpotent element $f\in\spa_{2n}$
in the form
\beql{fspn}
f=F_{2\tss 1}+F_{3\tss 2}+\dots+F_{n\tss n-1}+\frac12\ts F_{n'\tss n}.
\eeq
The $\sll_2$-triple is formed by the elements $\{e,f,h\}$ with
\ben
e=\sum_{i=1}^{n-1} i\tss(2\tss n-i)\ts F_{i\ts i+1}+\frac{n^2}2\ts F_{n\tss n'}\Fand
h=\sum_{i=1}^{n} (2\tss n-2\tss i+1)\ts F_{i\tss i}.
\een
The invariant symmetric bilinear form on $\spa_{2n}$
is defined by
\beql{invformc}
(X|Y)=\frac12\ts \tr\ts XY,\qquad X,Y\in\spa_{2n},
\eeq
where $X$ and $Y$ are understood as $2n\times 2n$ symplectic matrices over $\CC$.

Consider the determinant \eqref{det} of the matrix with entries in $\Vc(\p)\ot\CC[\di]$,
\ben
\det
\begin{bmatrix}\di+F_{1\tss 1}&1&\dots&0&0&0&\dots&0&0\\
                 F_{2\tss 1}&\di+F_{2\tss 2}&\dots&0&0&0&\dots&0&0\\
                 \dots&\dots&\ddots&\dots&\dots&\dots&\dots&\dots&\dots\\
                 F_{n\tss 1}&F_{n\tss 2}&\dots&\di+F_{n\tss n}&1&0&\dots&0&0\\
                 F_{n^{\tss\prime}\tss 1}&F_{n^{\tss\prime}\tss 2}&\dots
                 &F_{n^{\tss\prime}\ts n}&\di+F_{n^{\tss\prime}\tss n^{\tss\prime}}
                 &-1\phantom{-}&\dots&0&0 \\
                 \dots&\dots&\dots&\dots&\dots&\dots&\dots&\dots&\dots\\
                 F_{2^{\tss\prime}\tss 1}&F_{2^{\tss\prime}\tss 2}
                 &\dots&F_{2^{\tss\prime}\tss n}&F_{2^{\tss\prime}\tss n^{\tss\prime}}
                          &\dots&\dots&\di+F_{2^{\tss\prime}\tss 2^{\tss\prime}}&-1\phantom{-}\\
                 F_{1^{\tss\prime}\tss 1}&F_{1'\tss 2}&\dots&F_{1'\tss n}
                 &F_{1'\tss n^{\tss\prime}}
                           &\dots&\dots&F_{1'\tss 2^{\tss\prime}}&\di+F_{1'\tss 1'}
                \end{bmatrix}
\een
so that the $(i,j)$ entry of the matrix is $\de_{ij}\tss\di+F_{i\tss j}$ for $i\geqslant j$,
the $(i,i+1)$ entry is $1$ for $i\leqslant n$ and $-1$ for $i>n$, while the remaining entries
are zero. The determinant has the form
\beql{diffsp}
\di^{\ts 2n}+w_2\ts \di^{\ts 2n-2}+w_3\ts \di^{\ts 2n-3}+\dots+w_{2n},
\qquad w_i\in\Vc(\p).
\eeq

\bth\label{thm:c}
All elements $w_2,\dots,w_{2n}$ belong to the classical $\Wc$-algebra
$\Wc(\spa_{2n},f)$. Moreover, the elements $w^{(r)}_2,w^{(r)}_4,\dots,w^{(r)}_{2n}$
with $r=0,1,\dots$ are algebraically independent
and generate the algebra $\Wc(\spa_{2n},f)$.
\eth

\bpf
Denote the determinant by $D$ and let $D_i$ (resp., $\overline D_i$) denote
the $i\times i$ minor
corresponding to the first (resp., last) $i$ rows and columns. We suppose that
$D_0=\overline D_0=1$. Lemma~\ref{lem:deex} implies the expansion
\ben
D=D_n\ts \overline D_n+\sum_{j,k=1}^{n}(-1)^{n-j+1}\ts
D_{j-1}\ts F_{k'\tss j}\ts \overline D_{k-1}.
\een
If $1\leqslant i\leqslant n-1$ then the relation
$\rho\tss\{F_{i\ts i+1}\ts^{}_{\la}\ts D\}=0$
follows by the same calculation as in the proof of Theorem~\ref{thm:b}.
Furthermore,
\begin{multline}
\rho\tss\{F_{n\ts n'}\ts^{}_{\la}\ts D\}=-2\tss D^+_{n-1}\ts \overline D_{n}
+2\tss D^+_{n}\ts \overline D_{n-1}
+2\ts\sum_{k=1}^{n-1}D^+_{n-1}\ts F_{k'\tss n'}\ts \overline D_{k-1}\\
-2\ts D^+_{n-1}\ts (2\tss F_{n\tss n}+\la)\ts \overline D_{n-1}
+2\ts\sum_{j=1}^{n-1}(-1)^{n-j+1}\ts D^+_{j-1}\ts F_{n\ts j}\ts \overline D_{n-1}.
\non
\end{multline}
This is zero since the relations \eqref{dplusex} and \eqref{dbarplusex}
are valid for the case of $\spa_{2n}$ as well.
Thus,
all elements $w_2,\dots,w_{2n}$ belong to the subalgebra $\Wc(\spa_{2n},f)$
of $\Vc(\p)$.

The odd powers of the matrix $f$ form a basis of the centralizer
$\spa_{2n}^{\tss f}$.
For $j=1,\dots,n$ take $v_j$ to be equal, up to a sign, to $f^{2\tss j-1}$.
The condition of Proposition~\ref{prop:genwalg} will hold for the family
of elements $w_2,w_4,\dots,w_{2n}$,
thus implying that they are generators
of the differential algebra $\Wc(\spa_{2n},f)$.

The images of the elements $w_k$ under the homomorphism \eqref{phiw}
are the elements $\wt w_k\in\Vc(\h)$ found from the relation
\begin{multline}
(\di+F_{11})\dots (\di+F_{n\tss n})\ts
(\di+F_{n'\tss n'}) \dots (\di+F_{1'\tss 1'})\\[0.3em]
=\di^{\ts 2n}+\wt w_2\ts \di^{\ts 2n-2}+\wt w_3\ts \di^{\ts 2n-3}+\dots+\wt w_{2n}.
\non
\end{multline}
In the notation of Sec.~\ref{subsec:ctt} we have
\ben
h^{(r)}_j=F^{(r)}_{j\tss j}-F^{(r)}_{j+1\ts j+1},\quad j=1,\dots,n-1,\Fand
h^{(r)}_n=F^{(r)}_{n\tss n}.
\een
The Cartan matrix is of the size $n\times n$,
\ben
A=\begin{bmatrix}2&-1\phantom{-}&0&\dots&0&0\\
                 -1\phantom{-}&2&-1\phantom{-}&\dots&0&0\\
                 0&-1\phantom{-}&2&\dots&0&0\\
                 \dots&\dots&\dots&\ddots&\dots &\dots \\
                             0&0&0&\dots&2&-2\phantom{-}\\
                             0&0&0&\dots&-1\phantom{-}&2
                \end{bmatrix}
\een
so that $a_{ii}=2$ for $i=1,\dots,n$
and $a_{i\ts i+1}=a_{i+1\ts i}=a_{n\ts n-1}=-1$ for $i=1,\dots,n-2$, while
$a_{n-1\ts n}=-2$ and
all other
entries are zero. The entries of the diagonal matrix $D=\diag\tss[\ep_1,\dots,\ep_n]$
are found by
\ben
\ep_1=\dots=\ep_{n-1}=1\Fand \ep_n=1/2.
\een
Hence, regarding $\Vc(\h)$ as the algebra of polynomials
in the variables $F^{(r)}_{i\tss i}$ with $i=1,\dots,n$ and $r=0,1,\dots$, we get
\ben
\sum_{j=1}^{n} a_{ji}\ts\frac{\di}{\di\tss h^{(r)}_j}=
\frac{\di}{\di\tss F^{(r)}_{i\tss i}}-\frac{\di}{\di\tss F^{(r)}_{i+1\ts i+1}},
\qquad i=1,\dots,n-1,
\een
and
\ben
\sum_{j=1}^{n} a_{j\tss n}\ts\frac{\di}{\di\tss h^{(r)}_j}=
2\ts\frac{\di}{\di\tss F^{(r)}_{n\tss n}}.
\een
Therefore, the screening operators \eqref{scr} take the form
\ben
V_i=\sum_{r=0}^{\infty} V_{i\ts r}
\Bigg(\frac{\di}{\di\tss F^{(r)}_{i\tss i}}
-\frac{\di}{\di\tss F^{(r)}_{i+1\ts i+1}}\Bigg),\qquad i=1,\dots,n-1,
\een
and
\ben
V_n=2\ts\sum_{r=0}^{\infty} V_{n\ts r}\ts\frac{\di}{\di\tss F^{(r)}_{n\tss n}},
\een
where the coefficients $V_{i\ts r}$ are found by the
relations
\ben
\sum_{r=0}^{\infty}\frac{V_{i\ts r}\ts z^r}{r!}=\exp\Big({-}\sum_{m=1}^{\infty}
\frac{F^{(m-1)}_{i\tss i}-F^{(m-1)}_{i+1\ts i+1}}{m!}\ts z^m\Big),\qquad i=1,\dots,n-1,
\een
and
\ben
\sum_{r=0}^{\infty}\frac{V_{n\ts r}\ts z^r}{r!}=\exp\Big({-}\sum_{m=1}^{\infty}
\frac{2\ts F^{(m-1)}_{n\tss n}}{m!}\ts z^m\Big).
\een
The differential algebra $\wt\Wc(\spa_{2n},f)$ consists of the polynomials
in the variables $F^{(r)}_{i\tss i}$, which are annihilated by all
operators $V_i$. One easily verifies directly
that the elements $\wt w_k$ belong to $\wt\Wc(\spa_{2n},f)$;
cf. \cite[Sec.~4.2]{mm:yc}.
Moreover, the elements $\wt w^{\tss(r)}_2,\wt w^{\tss(r)}_4,\dots,
\wt w^{\tss(r)}_{2n}$
with $r$ running over nonnegative integers
are algebraically independent
generators of the algebra $\wt\Wc(\spa_{2n},f)$; see \cite[Ch.~8]{f:lc}.
Hence, the generators $w^{\tss(r)}_2,w^{\tss(r)}_4,\dots,
w^{\tss(r)}_{2n}$ of $\Wc(\spa_{2n},f)$
are also algebraically independent.
\epf

The injective homomorphism $\Wc(\spa_{2n},f)\hra \Vc(\h)$
taking $w_i$ to $\wt w_i$
is known as the {\it Miura transformation\/} in type $C$ \cite{ds:la};
see also \cite[Ch.~8]{f:lc}. Note that by the arguments used in
the proof of Theorem~\ref{thm:c},
the homomorphism \eqref{phiw}
is bijective in the case $\g=\spa_{2n}$.
An alternative proof of the
first part of the theorem can be obtained with the use of the folding
procedure; see Remark~\ref{rem:fold}.

Another family of generators
of the classical $\Wc$-algebra $\Wc(\spa_{2n},f)$ can be constructed
by using Proposition~\ref{prop:mmta}; cf. Corollaries~\ref{cor:genh}
and~\ref{cor:genhb}. Extend the definition \eqref{wtab}
to the symplectic case by restricting the range of the indices
to $2n\geqslant i\geqslant j\geqslant 1$.
Furthermore, define elements $e_{m\ts 0},h_{m\ts 0}\in \Vc(\p)$
by the formulas \eqref{emob} and \eqref{hmob}.
Proposition~\ref{prop:ewcoin} implies
$e_{m\ts 0}=w_{m}$ for $m=2,3,\dots,2n$.
Hence, by Theorem~\ref{thm:c} the family
$e^{(r)}_{2\ts 0},e^{(r)}_{4\ts 0},\dots,e^{(r)}_{2n\ts 0}$
with $r=0,1,\dots$
is algebraically independent
and generates the algebra $\Wc(\spa_{2n},f)$.

\bco\label{cor:genhc}
All elements $h_{m\ts 0}$ belong to the classical $\Wc$-algebra $\Wc(\spa_{2n},f)$.
Moreover, the family
$h^{(r)}_{2\ts 0},h^{(r)}_{4\ts 0},\dots,h^{(r)}_{2n\ts 0}$
with $r=0,1,\dots$
is algebraically independent
and generates $\Wc(\spa_{2n},f)$.
\qed
\eco

\section{Generators of $\Wc(\oa_{2n},f)$}
\label{sec:gend}
\setcounter{equation}{0}

We keep the notation for the generators of the Lie
algebra $\oa_{2n}$ introduced in the beginning of Sec.~\ref{sec:genb}
by taking $N=2n$. In particular, we have the relations \eqref{on}
and \eqref{symrel}. We will work with the algebra of
pseudo-differential operators $\Vc(\p)\ot\CC((\di^{-1}))$,
where the relations are given by \eqref{dfpol} and
\ben
\di^{-1} \tss F^{(r)}_{ij}=\sum_{s=0}^{\infty} (-1)^s\ts F^{(r+s)}_{ij}\ts \di^{-s-1}.
\een

Take the principal nilpotent element $f\in\oa_{2n}$
in the form
\beql{fdn}
f=F_{2\tss 1}+F_{3\tss 2}+\dots+F_{n\tss n-1}+F_{n'\tss n-1}.
\eeq
The $\sll_2$-triple is formed by the elements $\{e,f,h\}$ with
\ben
e=\sum_{i=1}^{n-2} i\tss(2\tss n-i-1)\ts F_{i\ts i+1}
+\frac{n^2-n}{2}\ts\big(F_{n-1\tss n}+F_{n-1\tss n'}\big)\fand
h=2\ts\sum_{i=1}^{n-1} (n-i)\ts F_{i\tss i}.
\een
The invariant symmetric bilinear form on $\oa_{2n}$
is defined by
\beql{invformd}
(X|Y)=\frac12\ts \tr\ts XY,\qquad X,Y\in\oa_{2n},
\eeq
where $X$ and $Y$ are understood as matrices over $\CC$
which are skew-symmetric with respect to the antidiagonal.

Consider the following $(2n+1)\times(2n+1)$ matrix with entries in
$\Vc(\p)\ot\CC((\di^{-1}))$,
\ben
\begin{bmatrix}\di+F_{1\tss 1}&1&\dots&0&0&0&\dots&0&0\\
                 F_{2\tss 1}&\di+F_{2\tss 2}&\dots&0&0&0&\dots&0&0\\
                 \dots&\dots&\ddots&\dots&\dots&\dots&\dots&\dots&\dots\\
                 F_{n\tss 1}-F_{n'\tss 1}&F_{n\tss 2}-F_{n'\tss 2}&\dots&
                 \di+F_{n\tss n}&0&-2\tss\di\phantom{-}&\dots&0&0\\
                 0&0&\dots&0&\phantom{-}\di^{-1}&0&\dots&0&0 \\
                 F_{n^{\tss\prime}\tss 1}&F_{n^{\tss\prime}\tss 2}
                            &\dots&0&0
                            &\di+F_{n^{\tss\prime}\tss n^{\tss\prime}}&\dots&0&0\\
                 \dots&\dots&\dots&\dots&\dots&\dots&\ddots&\dots&\dots\\
                 F_{2^{\tss\prime}\tss 1}&0&\dots&\dots&0
                          &F_{2'\tss n'}-F_{2'\tss n}
                          &\dots&\di+F_{2^{\tss\prime}\tss 2^{\tss\prime}}&-1\phantom{-}\\
                 0&F_{1'\tss 2}&\dots&\dots&0
                           &F_{1'\tss n'}-F_{1'\tss n}
                           &\dots&F_{1'\tss 2^{\tss\prime}}&\di+F_{1'\tss 1'}
                \end{bmatrix}
\een
so that the first $n-1$ rows, the last $n-1$ columns and the lower-left $n\times n$
submatrix respectively
coincide with those of the matrix introduced for $\oa_{2n+1}$ in Sec.~\ref{sec:genb}.
All entries in the row and column $n+1$ are zero, except for the $(n+1,n+1)$ entry
which equals $\di^{-1}$. The
$(n,j)$ entries are
$F_{n\tss j}-F_{n'\tss j}$ for $j=1,\dots,n-1$, the $(n,n)$ entry is $\di+F_{n\tss n}$
and the $(n,n+2)$ entry is $-2\tss\di$, while the remaining entries in row $n$ are zero.
Finally, the remaining nonzero entries in column $n+2$ are $F_{k'\tss n'}-F_{k'\tss n}$
for $k=1,2,\dots,n-1$ which occur in the respective rows $2\tss n-k+2$, and
$\di+F_{n'\tss n'}$ which occurs in row $n+2$.

One easily verifies that the column-determinant and row-determinant of
this matrix coincide, so that the determinant \eqref{det}
is well-defined and we denote it by $D$.
Applying the simultaneous column expansion along the first $n$ columns and
using Lemma~\ref{lem:deex}, we derive that
it can be written in the form
\beql{diffd}
D=D_n\ts \di^{-1}\ts\overline D_n+2\ts\sum_{j,k=1}^{n}(-1)^{n-j}\ts
D_{j-1}\ts F_{k'\tss j}\ts \overline D_{k-1},
\eeq
where $D_i$ (resp., $\overline D_i$) denotes
the $i\times i$ minor
corresponding to the first (resp., last) $i$ rows and columns. We suppose that
$D_0=\overline D_0=1$. Write
\ben
\bal
D_n&=\di^{\ts n}+y_1\ts \di^{\ts n-1}+y_2\ts \di^{\ts n-2}+\dots+y_{n},\\[0.3em]
\overline D_n&=\di^{\ts n}+\di^{\ts n-1}\ts \bar y_1+
\di^{\ts n-2}\ts \bar y_2+\dots+\bar y_{n},
\eal
\een
for certain uniquely determined elements $y_i,\bar y_i\in\Vc(\p)$.

\ble\label{lem:ddbar}
We have $\bar y_i=(-1)^i\ts y_i$ for all $i=1,\dots,n$.
\ele

\bpf
Replace $\di$ by $-\di$ in the minor $\overline D_n$
and multiply each row by $-1$. The lemma can then be equivalently stated as
the identity
\ben
(-1)^n\ts \overline D_n\Big|_{\di\ts\mapsto\ts -\di}=\di^{\ts n}+\di^{\ts n-1}\ts y_1+
\di^{\ts n-2}\ts y_2+\dots+y_{n}.
\een
The left hand side is the determinant
\ben
\det
\begin{bmatrix}\di+F_{n\tss n}&1&0&0&\dots&0\\
                 F_{n\ts n-1}-F_{n'\ts n-1}&\di+F_{n-1\ts n-1}&1&0&\dots&0\\
                 \dots&\dots&\dots&\dots&\dots& \dots \\
                 F_{n\ts 2}-F_{n'\ts 2}&F_{n-1\tss 2}&\dots&\dots&\di+F_{2\tss 2}&1\\
                 F_{n\ts 1}-F_{n'\ts 1}&F_{n-1\tss 1}&\dots&\dots&F_{2\tss 1}&\di+F_{1\tss 1}
                \end{bmatrix}.
\een
The claim now follows from the observation that this determinant coincides with
the image of $D_n$ under the anti-automorphism of the algebra $\Vc(\p)\ot \CC[\di]$
which is the identity on the generators $F_{ij}$ and $\di$.
\epf

Lemma~\ref{lem:ddbar} implies that
the
pseudo-differential operator $D$ can be written as
\beql{diffod}
D=\di^{\ts 2n-1}+w_2\ts \di^{\ts 2n-3}+w_3\ts \di^{\ts 2n-4}+\dots+w_{2n-1}
+(-1)^n\ts y_n\ts \di^{-1}\ts y_n,\qquad w_i\in\Vc(\p).
\eeq

\bth\label{thm:d}
The elements $w_2,w_3,\dots,w_{2n-1}$ and $y_n$
belong to the classical $\Wc$-algebra
$\Wc(\oa_{2n},f)$. Moreover,
the elements $w^{(r)}_2,w^{(r)}_4,\dots,w^{(r)}_{2n-2},y^{(r)}_n$
with $r=0,1,\dots$ are algebraically independent
and generate the algebra $\Wc(\oa_{2n},f)$.
\eth

\bpf
The relation
$\rho\tss\{F_{i\ts i+1}\ts^{}_{\la}\ts D\}=0$ for
$1\leqslant i\leqslant n-1$
follows by the same calculations as in the proof of Theorem~\ref{thm:b}.
Furthermore, let $\si=(n\ts n')$ be the permutation of the set of indices $\{1,\dots,2n\}$
which swaps $n$ and $n'=n+1$ and leaves all other indices fixed.
The mapping
\beql{invol}
\varsigma:F_{i\tss j}\mapsto F_{\si(i)\tss \si(j)}
\eeq
defines an involutive automorphism of the Lie algebra $\oa_{2n}$.
It also extends to an involutive automorphism of the Poisson vertex algebra
$\Vc(\oa_{2n})$. We claim that all coefficients of the pseudo-differential operator $D$
are $\varsigma$-invariant. Indeed, let us apply the following operations
on rows and columns of the given matrix. Replace row $n+2$ by the sum
of rows $n$ and $n+2$. Then replace column $n$ by the sum
of columns $n$ and $n+2$. Finally, multiply row $n$ and column $n+2$ by $-1$.
As a result, we get the image of the matrix with respect to the involution \eqref{invol}.
On the other hand, the determinant $D$ remains unchanged. This proves the relation
$\rho\tss\{F_{n-1\ts n'}\ts^{}_{\la}\ts D\}=0$. This shows that all coefficients
of the operator $D$ belong to $\Wc(\oa_{2n},f)$.

Now consider the minor $D_n$. It can be written in the form
\begin{multline}
D_n=D_{n-1}\ts (\di+F_{n\tss n})-D_{n-2}\ts (F_{n\tss n-1}-F_{n'\tss n-1})
+D_{n-3}\ts (F_{n\tss n-2}-F_{n'\tss n-2})\\[0.3em]
{}+\dots+(-1)^{n-2} \ts D_1\tss (F_{n\tss 2}-F_{n'\tss 2})+(-1)^{n-1}
\ts D_0\tss (F_{n\tss 1}-F_{n'\tss 1}).
\label{dnexpd}
\end{multline}
The same calculations as in the proof of Theorem~\ref{thm:gln},
show that $\rho\tss\{F_{i\ts i+1}\ts^{}_{\la}\ts D_n\}=0$ for $i=1,\dots,n-1$.
Furthermore,
\begin{multline}
\rho\tss\{F_{n-1\ts n'}\ts^{}_{\la} D_n\}=-D^+_{n-1}
-D^+_{n-2}\ts  (\di+F_{n\tss n})+D^+_{n-2}\ts  (F_{n-1\ts n-1}+F_{n\tss n}+\la)\\[0.3em]
-D^+_{n-3}\ts F_{n-1\ts n-2}
{}+\dots+(-1)^{n-1} \ts D^+_1\tss F_{n-1\ts 2}+(-1)^{n} \ts D^+_0\tss F_{n-1\tss 1}.
\non
\end{multline}
Applying relation \eqref{dnexpa} to the determinant $D^+_{n-1}$
we get
\ben
\rho\tss\{F_{n-1\ts n'}\ts^{}_{\la} D_n\}=-2\ts D^+_{n-2}\ts \di.
\een
This implies $\rho\tss\{F_{n-1\ts n'}\ts^{}_{\la}\ts y_n\}=0$ so that
the constant term $y_n$ of the differential operator
$D_n$ belongs to $\Wc(\oa_{2n},f)$.

To use Proposition~\ref{prop:genwalg}, note that
the odd powers $f, f^3,\dots,f^{2n-3}$ of the matrix $f$ together
with the element $F_{n\tss 1}-F_{n'\tss 1}$
form a basis of the centralizer
$\oa_{2n}^{\tss f}$. Take $v_j$ to be equal, up to a sign, to $f^{2\tss j-1}$
for $j=1,\dots,n-1$ and set $v_n=F_{n\tss 1}-F_{n'\tss 1}$.
The condition of Proposition~\ref{prop:genwalg} holds for the family
$w_2,w_4,\dots,w_{2n-2},y_n$,
thus implying that they are generators
of the differential algebra $\Wc(\oa_{2n},f)$.

The images of the coefficients of the operator $D$
under the homomorphism \eqref{phiw}
are the elements $\wt w_k,\wt y_n\in\Vc(\h)$ found from the relation
\begin{multline}
(\di+F_{11})\dots (\di+F_{n\tss n})\ts\di^{-1}\ts
(\di+F_{n'\tss n'}) \dots (\di+F_{1'\tss 1'})\\[0.3em]
=\di^{\ts 2n-1}+\wt w_2\ts \di^{\ts 2n-3}+\wt w_3\ts \di^{\ts 2n-4}+\dots+
\wt w_{2n-1}+(-1)^n\ts \wt y_n\ts \di^{-1}\ts \wt y_n.
\non
\end{multline}
In particular,
\ben
\wt y_n=(\di+F_{11})\dots (\di+F_{n\tss n})\ts 1.
\een

In the notation of Sec.~\ref{subsec:ctt} we have
\ben
h^{(r)}_j=F^{(r)}_{j\tss j}-F^{(r)}_{j+1\ts j+1},\quad j=1,\dots,n-1,\Fand
h^{(r)}_n=F^{(r)}_{n-1\ts n-1}+F^{(r)}_{n\tss n}.
\een
The Cartan matrix is of the size $n\times n$,
\ben
A=\begin{bmatrix}2&-1\phantom{-}&0&\dots&0&0&0\\
                 -1\phantom{-}&2&-1\phantom{-}&\dots&0&0&0\\
                 0&-1\phantom{-}&2&\dots&0&0&0\\
                 \dots&\dots&\dots&\ddots&\dots&\dots &\dots \\
                 0&0&0&\dots&2&-1\phantom{-}&-1\phantom{-}\\
                             0&0&0&\dots&-1\phantom{-}&2&0\\
                             0&0&0&\dots&-1\phantom{-}&0&2
                \end{bmatrix}
\een
so that $a_{ii}=2$ for $i=1,\dots,n$
and $a_{i\ts i+1}=a_{i+1\ts i}=a_{n-2\ts n}=a_{n\ts n-2}=-1$ for $i=1,\dots,n-2$
and all other
entries are zero. Then the diagonal matrix $D$ (see Sec.~\ref{subsec:ctt})
is the identity matrix so that
$
\ep_1=\dots=\ep_{n}=1.
$
Hence, regarding $\Vc(\h)$ as the algebra of polynomials
in the variables $F^{(r)}_{i\tss i}$ with $i=1,\dots,n$ and $r=0,1,\dots$, we get
\ben
\sum_{j=1}^{n} a_{ji}\ts\frac{\di}{\di\tss h^{(r)}_j}=
\frac{\di}{\di\tss F^{(r)}_{i\tss i}}-\frac{\di}{\di\tss F^{(r)}_{i+1\ts i+1}},
\qquad i=1,\dots,n-1,
\een
and
\ben
\sum_{j=1}^{n} a_{j\tss n}\ts\frac{\di}{\di\tss h^{(r)}_j}=
\frac{\di}{\di\tss F^{(r)}_{n-1\ts n-1}}+\frac{\di}{\di\tss F^{(r)}_{n\tss n}}.
\een
Therefore, the screening operators \eqref{scr} take the form
\ben
V_i=\sum_{r=0}^{\infty} V_{i\ts r}
\Bigg(\frac{\di}{\di\tss F^{(r)}_{i\tss i}}
-\frac{\di}{\di\tss F^{(r)}_{i+1\ts i+1}}\Bigg),\qquad i=1,\dots,n-1,
\een
and
\ben
V_n=\sum_{r=0}^{\infty} V_{n\ts r}\ts
\Bigg(\frac{\di}{\di\tss F^{(r)}_{n-1\ts n-1}}
+\frac{\di}{\di\tss F^{(r)}_{n\tss n}}\Bigg),
\een
where the coefficients $V_{i\ts r}$ are found by the
relations
\ben
\sum_{r=0}^{\infty}\frac{V_{i\ts r}\ts z^r}{r!}=\exp\Big({-}\sum_{m=1}^{\infty}
\frac{F^{(m-1)}_{i\tss i}-F^{(m-1)}_{i+1\ts i+1}}{m!}\ts z^m\Big),\qquad i=1,\dots,n-1,
\een
and
\ben
\sum_{r=0}^{\infty}\frac{V_{n\ts r}\ts z^r}{r!}=\exp\Big({-}\sum_{m=1}^{\infty}
\frac{F^{(m-1)}_{n-1\ts n-1}+F^{(m-1)}_{n\tss n}}{m!}\ts z^m\Big).
\een
The differential algebra $\wt\Wc(\oa_{2n},f)$ consists of the polynomials
in the variables $F^{(r)}_{i\tss i}$, which are annihilated by all
operators $V_i$. It is easy to verify directly
that the elements $\wt w_k$ and $\wt y_n$ belong to $\wt\Wc(\oa_{2n},f)$;
cf. \cite[Sec.~4.2]{mm:yc}.
The elements $\wt w^{\tss(r)}_2,\wt w^{\tss(r)}_4,\dots,
\wt w^{\tss(r)}_{2n-2},\wt y^{\tss(r)}_{n}$
with $r$ running over nonnegative integers
are known to be algebraically independent
generators of the algebra $\wt\Wc(\oa_{2n},f)$; see \cite[Ch.~8]{f:lc}.
Hence, the generators $w^{\tss(r)}_2,w^{\tss(r)}_4,\dots,
w^{\tss(r)}_{2n-2},y^{\tss(r)}_{n}$ of $\Wc(\oa_{2n},f)$
are also algebraically independent.
\epf

By the arguments used in
the proof of Theorem~\ref{thm:d},
the homomorphism \eqref{phiw}
is bijective in the case $\g=\oa_{2n}$.

Another family of generators
of $\Wc(\oa_{2n},f)$ analogous to those described in Corollaries~\ref{cor:genh},
\ref{cor:genhb} and \ref{cor:genhc}, can be
obtained from an appropriate analogue of Proposition~\ref{prop:mmta}.
This is yet another version of the MacMahon Master Theorem
involving the noncommutative elementary and complete symmetric functions
associated with the $(2n+1)\times(2n+1)$ generators matrix introduced
in the beginning of this section.

Using the determinant $D=D(\di)$ define the elements $e_m\in\Vc(\p)\ot\CC[\di]$
as the coefficients of the formal power series in $t$,
\ben
\sum_{m=0}^{\infty}e_m\ts t^m=t^{\tss 2\tss n-1} D(\di+t^{-1}).
\een
Denote this series by $e(t)$. Furthermore, define the elements
$h_m\in\Vc(\p)\ot\CC[\di]$ as the coefficients of the series
\ben
h(t)=\sum_{m=0}^{\infty}h_m\ts t^m, \qquad h(t)=e(-t)^{-1}.
\een
Write the expansions \eqref{emexp} and \eqref{hmexp}.
By Theorem~\ref{thm:d}, all coefficients of the
differential operators $e_m$ and $h_m$ belong to $\Wc(\oa_{2n},f)$.
Explicit formulas for these coefficients
in terms of the generators $F_{ij}$ and $\di$
take the form similar to those of Sec.~\ref{subsec:mmt}
but are more complicated. For this reason we will not reproduce them.
On the other hand, their images under the isomorphism \eqref{phiw} are easy
to describe. Set $a_{i\tss i}=\di+F_{i\tss i}$. We have
\ben
\phi:e(t)\mapsto (1+t\ts a_{1\tss 1})\dots
(1+t\ts a_{n\tss n})\ts(1+t\ts \di)^{-1}\ts (1+t\ts a_{n'\tss n'})
\dots (1+t\ts a_{1'\tss 1'})
\een
and hence
\ben
\phi:h(t)\mapsto (1-t\ts a_{1'\tss 1'})^{-1}\dots
(1-t\ts a_{n'\tss n'})^{-1}\ts(1-t\ts \di)\ts (1-t\ts a_{n\tss n})^{-1}
\dots (1-t\ts a_{1\tss 1})^{-1}.
\een
In particular, the constant terms $e_{m\tss 0}$ of the differential
operators $e_m$ coincide with the elements $w_m$ for $m=2,\dots,2n-1$.
Observe that since $a_{n\tss n}+a_{n'\tss n'}=2\tss \di$, we have the relation
\ben
(1-t\ts a_{n'\tss n'})^{-1}\ts(1-t\ts \di)\ts (1-t\ts a_{n\tss n})^{-1}
=\frac12\ts \Big((1-t\ts a_{n\tss n})^{-1}+(1-t\ts a_{n'\tss n'})^{-1}\Big).
\een
Therefore, the images of the elements $h_m$ can be written explicitly as
\ben
\phi(h_m)=\frac12\ts \sum_{k_{1'}+\dots+k_1=m}
a_{1'\tss 1'}^{k_{1'}}\dots a_{n'\tss n'}^{k_{n'}}
\ts a_{n\tss n}^{k_{n}} \dots a_{1\tss 1}^{k_{1}},
\een
summed over nonnegative integers $k_i$ such that either $k_{n}=0$ or $k_{n'}=0$.

The generators $\phi(w_i)$ and $\phi(y_n)$ of the differential algebra $\wt\Wc(\oa_{2n},f)$
were introduced in \cite{ds:la} in relation with the {\it Miura transformation\/}
in type $D$. The coefficients of the differential operators
$\phi(h_m)$ were found in \cite{mm:yc}
as the Harish-Chandra images of generators of the center of the affine vertex
algebra at the critical level in type $D$ (in a slightly different notation).
The relation $e(-t)\tss h(t)=1$ thus makes a connection between the two
families of generators. Moreover, we obtain the following.

\bco\label{cor:genhd}
The family
$h^{(r)}_{2\ts 0},h^{(r)}_{4\ts 0},\dots,h^{(r)}_{2n-2\ts 0},y^{(r)}_n$
with $r=0,1,\dots$
is algebraically independent
and generates $\Wc(\oa_{2n},f)$.
\qed
\eco

\section{Generators of $\Wc(\g_{2},f)$}
\label{sec:geng}
\setcounter{equation}{0}

We denote by $\g_2$ the simple Lie algebra of type $G_2$ with the Cartan matrix
\ben
A=\begin{bmatrix}\phantom{-}2&-1\phantom{-}\\
                -3&2
                \end{bmatrix}.
\een
This Lie algebra is $14$-dimensional. An explicit basis and
the multiplication table is given in \cite[Lec.~22]{fh:rt}.
A simple combinatorial construction of $\g_2$ is given
in \cite{w:ccg} which we will follow below. We let $\al$ and $\be$ denote
the simple roots, and the set of positive roots is
\ben
\al,\quad\be,\quad\al+\be,\quad\al+2\tss \be,\quad\al+3\tss \be,\quad 2\tss \al+3\tss \be.
\een
For each positive root $\ga$ we let $X_{\ga}$ and $Y_{\ga}$ denote the root
vectors associated with $\ga$ and $-\ga$, respectively. The Cartan subalgebra
$\h$ of $\g_2$ is two-dimensional and we let $H_{\al}$ and $H_{\be}$ denote
its basis elements. The full list of commutators between the basis elements
of $\g_2$ is given in \cite{w:ccg}.\footnote{The preprint version of \cite{w:ccg}
provides such relations for the basis elements which differ from the journal
version by signs.} In particular,
\ben
[X_{\al}, Y_{\al}]=H_{\al},\qquad [X_{\be}, Y_{\be}]=H_{\be}
\een
and
\ben
\bal[]
[H_{\al}, X_{\al}]&=-2\tss X_{\al},\qquad [H_{\al}, X_{\be}]=X_{\be},\\
[H_{\be}, X_{\be}]&=-2\tss X_{\be},\qquad [H_{\be}, X_{\al}]=3\tss X_{\al}.
\eal
\een
The commutation relations can also be recovered from the identification
of the basis elements of $\g_2$ as elements of $\oa_7$ or $\oa_8$
under the embeddings $\g_2\subset \oa_7$ or $\g_2\subset \oa_8$ which
we use below.

The respective subsets of elements $X_{\ga}$ and $Y_{\ga}$
with $\ga$ running over the set of positive roots
span the nilpotent subalgebras $\n_+$ and $\n_-$.
As before, we set $\p=\n_-\oplus\h$.
We will be working with the algebra of differential operators $\Vc(\p)\ot\CC[\di]$,
where
\ben
\di\ts X^{(r)}- X^{(r)}\ts \di=X^{\tss(r+1)}, \qquad X\in\p.
\een
Take the principal nilpotent element $f\in\g_{2}$
in the form
\beql{fgn}
f=Y_{\al}+Y_{\be}.
\eeq
The $\sll_2$-triple is formed by the elements $\{e,f,h\}$ with
\ben
e=-10\tss X_{\al}-6\tss X_{\be}\Fand
h=-10\tss H_{\al}-6\tss H_{\be}.
\een
The product $D^{-1}A$ is a symmetric matrix for $D=\diag[1,3]$ so
that $\ep_1=1$ and $\ep_2=3$.
The invariant symmetric bilinear form on $\g_2$
is then uniquely determined by the conditions
\beql{invformg}
(X_{\al}|\ts Y_{\al})=-1\Fand (X_{\be}|\ts Y_{\be})=-3.
\eeq

Consider the determinant \eqref{det} of the
$7\times 7$ matrix with entries in $\Vc(\p)\ot\CC[\di]$,
\ben
\det
\begin{bmatrix}\di+\wt F_{1\tss 1}&1&0&0&0&0&0\\[0.5em]
                 \frac13\ts Y_{\be}&\di+\wt F_{2\tss 2}&1&0&0&0&0\\[0.5em]
                 \frac13\ts Y_{\al+\be}&Y_{\al}&\di+\wt F_{3\tss 3}&1&0&0&0\\[0.5em]
                 \frac49\ts Y_{\al+2\tss\be}&-\frac23\ts Y_{\al+\be}
                 &\frac23\ts Y_{\be}&\di&-1\phantom{-}&0&0\\[0.5em]
                 -\frac49\ts Y_{\al+3\tss\be}&\frac49\ts Y_{\al+2\tss\be}
                 &0&-\frac23\ts Y_{\be}&\di-\wt F_{3\tss 3}&-1\phantom{-}&0\\[0.5em]
                 \frac49\ts Y_{2\tss\al+3\tss\be}&0
                            &-\frac49\ts Y_{\al+2\tss\be}
                            &\frac23\ts Y_{\al+\be}&-Y_{\al}
                            &\di-\wt F_{2\tss 2}&-1\phantom{-}\\[0.5em]
                 0&-\frac49\ts Y_{2\tss\al+3\tss\be}
                 &\frac49\ts Y_{\al+3\tss\be}&-\frac49\ts Y_{\al+2\tss\be}
                 &-\frac13\ts Y_{\al+\be}&-\frac13\ts Y_{\be}&\di-\wt F_{1\tss 1}
                \end{bmatrix},
\een
where
\ben
\wt F_{1\tss 1}=-H_{\al}-\frac23\ts H_{\be},\qquad \wt F_{2\tss 2}=-H_{\al}-\frac13\ts H_{\be}
\Fand \wt F_{3\tss 3}=-\frac13\ts H_{\be}.
\een
The determinant has the form
\beql{diffg}
D=\di^{\ts 7}+w_2\ts \di^{\ts 5}+w_3\ts \di^{\ts 4}+w_4\ts \di^{\ts 3}
+w_5\ts \di^{\ts 2}+w_6\ts \di+w_7,
\qquad w_i\in\Vc(\p).
\eeq

\bth\label{thm:g}
The elements $w_2,\dots,w_7$ belong to the classical $\Wc$-algebra
$\Wc(\g_{2},f)$. Moreover, the elements $w^{(r)}_2,w^{(r)}_6$
with $r=0,1,\dots$ are algebraically independent
and generate the algebra $\Wc(\g_{2},f)$.
\eth

\bpf
Consider the Lie algebra $\oa_8$ and its diagram automorphism of order $3$.
We can regard $\g_2$ as its fixed point subalgebra.
We will use a folding procedure to construct the classical $\Wc$-algebra
$\Wc(\g_{2},f)$ from its counterpart $\Wc(\oa_{8},f)$; cf. \cite{frs:fw}.
More precisely, using the notation of Sec.~\ref{sec:gend} for $\oa_8$,
identify $\g_2$ with a subalgebra of $\oa_8$ by setting
\ben
X_{\al}=-F_{2\tss 3},\qquad X_{\be}
=-F_{1\tss 2}-F_{3\tss 4}-F_{3\tss 4'},
\qquad Y_{\al}=F_{3\tss 2},
\qquad Y_{\be}=F_{2\tss 1}+F_{4\tss 3}+F_{4'\tss 3},
\een
where $i^{\tss\prime}=9-i$.
The remaining basis elements of $\g_2$ are then uniquely recovered,
and for the positive root vectors we have
\begin{alignat}{2}
X_{\al+\be}&=-F_{1\tss 3}+F_{2\tss 4}+F_{2\tss 4\pr},\qquad
&&X_{\al+2\tss\be}=-F_{1\tss 4}-F_{1\tss 4\pr}-F_{2\tss 3\pr},
\non\\
\qquad X_{\al+3\tss\be}&=F_{1\tss 3\pr},\qquad &&X_{2\tss\al+3\tss\be}=-F_{1\tss 2\pr}.
\non
\end{alignat}
The mapping $\psi:F_{i\tss j}\mapsto -F_{j\tss i}$ defines an involutive automorphism of
$\oa_8$. It preserves the subalgebra $\g_2$ and the
restriction to this subalgebra yields an involutive automorphism of the latter
such that
\ben
\psi:X_{\ga}\mapsto Y_{\ga},\qquad Y_{\ga}\mapsto X_{\ga},
\een
for all positive roots $\ga$. Moreover,
$\psi(H_{\al})=-H_{\al}$ and $\psi(H_{\be})=-H_{\be}$. So we thus obtain
expressions of the remaining basis vectors of $\g_2$ in terms of the generators
of $\oa_8$,
\ben
H_{\al}=-F_{2\tss 2}+F_{3\tss 3},\qquad H_{\be}=-F_{1\tss 1}+F_{2\tss 2}-2\tss F_{3\tss 3}
\een
and
\begin{alignat}{2}
Y_{\al+\be}&=F_{3\tss 1}-F_{4\tss 2}-F_{4\pr \tss 2},\qquad
&&Y_{\al+2\tss\be}=F_{4\tss 1}+F_{4\pr \tss 1}+F_{3\pr \tss 2},
\non\\
\qquad Y_{\al+3\tss\be}&=-F_{3\pr \tss 1},\qquad &&Y_{2\tss\al+3\tss\be}=F_{2\pr\tss 1}.
\non
\end{alignat}

Observe that the subalgebra $\g_2\subset \oa_8$ is contained in
the subalgebra $\oa_7\subset \oa_8$ defined as the span of the elements
\ben
F^{\circ}_{i\tss j}=F_{i\tss j},\qquad F^{\circ}_{i\tss j\pr}=\frac12\ts
F_{i\tss j\pr},\qquad F^{\circ}_{i\pr\tss j}=2\ts
F_{i\pr\tss j},\qquad 1\leqslant i,j\leqslant 3,
\een
and
\ben
F^{\circ}_{4\tss i}=F_{4\tss i}+F_{4'\tss i},
\qquad
F^{\circ}_{i\tss 4}=\frac12\ts\big(F_{i\tss 4}+F_{i\tss 4'}\big),\qquad
1\leqslant i\leqslant 3.
\een
In particular, we can identify $f\in\g_2$ with the principal nilpotent
elements
\ben
f=F^{\circ}_{2\tss 1}+F^{\circ}_{3\tss 2}+F^{\circ}_{4\tss 3}\in\oa_7
\Fand
f=F_{2\tss 1}+F_{3\tss 2}+F_{4\tss 3}+F_{4'\tss 3}\in\oa_8.
\een

Now we apply the folding procedure as a reduction from the classical
$\Wc$-algebra $\Wc(\oa_8,f)$ to $\Wc(\g_2,f)$; cf. Remark~\ref{rem:fold}.
Let $\vt$ denote the diagram automorphism of
the differential algebra $\Vc=\Vc(\oa_{8})$
defined by
\ben
\vt: F^{(r)}_{2\tss 3}\mapsto F^{(r)}_{2\tss 3},\quad
F^{(r)}_{1\tss 2}\mapsto F^{(r)}_{3\tss 4},\quad
F^{(r)}_{3\tss 4}\mapsto F^{(r)}_{3\tss 4'},\quad
F^{(r)}_{3\tss 4'}\mapsto F^{(r)}_{1\tss 2}
\een
and
\ben
\vt: F^{(r)}_{3\tss 2}\mapsto F^{(r)}_{3\tss 2},\quad
F^{(r)}_{2\tss 1}\mapsto F^{(r)}_{4\tss 3},\quad
F^{(r)}_{4\tss 3}\mapsto F^{(r)}_{4'\tss 3},\quad
F^{(r)}_{4'\tss 3}\mapsto F^{(r)}_{2\tss 1}.
\een
Note that $\vt^3=\text{\rm id}$ and we have
the direct sum decomposition
\ben
\oa_{8}=\g_2\oplus \oa^{(1)}_{8}\oplus \oa^{(2)}_{8},
\een
where $\g_2=\oa^{(0)}_{8}$ is
identified with the fixed point subalgebra under $\vt$, while
\ben
\oa^{(k)}_{8}=\{X\in\oa_{8}\ |\ \vt(X)=\om^k\ts X\},\qquad k=0,1,2,\qquad \om=e^{2\pi i/3},
\een
is the eigenspace of $\vt$ corresponding to the eigenvalue $\om^k$.
As in Sec.~\ref{subsec:ctt}, we will regard $\Vc=\Vc(\oa_8)$ as the symmetric
algebra $\Sr\big(t^{-1}\oa_{8}[t^{-1}]\big)$.
For the subalgebra of $\vt$-invariants in $\Vc$
we have
\ben
\Vc^{\ts\vt}=\Sr\Big(t^{-1}\g_2[t^{-1}]\oplus
\Sr^3\big(t^{-1}\oa^{(1)}_{8}[t^{-1}]\big)\oplus
\Sr^3\big(t^{-1}\oa^{(2)}_{8}[t^{-1}]\big)\oplus
t^{-1}\oa^{(1)}_{8}[t^{-1}]\ts t^{-1}\oa^{(2)}_{8}[t^{-1}]\Big).
\een
The $\la$-bracket on $\Vc$ defines
a $\la$-bracket on $\Vc^{\ts\vt}$ by restriction.
Let $\Jr$ be the ideal of $\Vc^{\ts\vt}$ generated by the
subspace
\ben
\Sr^3\big(t^{-1}\oa^{(1)}_{8}[t^{-1}]\big)\oplus
\Sr^3\big(t^{-1}\oa^{(2)}_{8}[t^{-1}]\big)\oplus
t^{-1}\oa^{(1)}_{8}[t^{-1}]\ts t^{-1}\oa^{(2)}_{8}[t^{-1}].
\een
For any element $P\in \Vc^{\ts\vt}$
we have $\{P^{}_{\ts\la\ts}\Jr\}\subset \Jr$.
The quotient space $\Vc^{\ts\vt}/\Jr$ is naturally identified
with the differential algebra $\Vc(\g_2)$.
This quotient is equipped with a $\la$-bracket induced from
that of $\Vc^{\ts\vt}$. The resulting bracket
on $\Vc(\g_2)$ is then obtained as the {\it folding\/}
of the $\la$-bracket on $\Vc(\oa_8)$.
It coincides with the $\la$-bracket
defined in \eqref{defla} for $\g=\g_2$.

For any element $P\in\Vc$ denote by $\text{\rm pr\tss}(P)$ its projection
to the subspace $\Vc^{\ts\vt}$
in the decomposition
\ben
\Vc=\Vc^{\ts\vt}\oplus t^{-1}\oa^{(1)}_{8}[t^{-1}]\ts\Vc^{\ts\vt}
\oplus t^{-1}\oa^{(2)}_{8}[t^{-1}]\ts\Vc^{\ts\vt}.
\een
Furthermore, let $\wt P\in\Vc(\g_2)$ be the image of $\text{\rm pr\tss}(P)$
under the natural
epimorphism
\ben
\Vc^{\ts\vt}\to \Vc^{\ts\vt}/\Jr\cong \Vc(\g_2).
\een
Now, if $P$ belongs to the classical $\Wc$-algebra $\Wc(\oa_{8},f)$,
then its image $\wt P$
belongs to the classical $\Wc$-algebra $\Wc(\g_2,f)$.
To prove the first part of the theorem, we will now verify
that the coefficients of the differential operator $D$ in \eqref{diffg}
are obtained as the respective images of the coefficients
of the pseudo-differential operator \eqref{diffd}.
Note that if $P\in\oa_8$, then we can write
\ben
3\tss P=P^{(0)}+P^{(1)}+P^{(2)},
\een
where
\ben
\bal
P^{(0)}&=P+\vt(P)+\vt^2(P)\in\g_2,\\
P^{(1)}&=P+\om^2\tss\vt(P)+\om\tss\vt^2(P)\in\oa^{(1)}_{8},\\
P^{(2)}&=P+\om\tss\vt(P)+\om^{2}\tss\vt^2(P)\in\oa^{(2)}_{8}.
\eal
\een
Hence, for the image $\wt P$ we have $3\tss \wt P=P^{(0)}$.
This gives the following formulas for the images of the generators of $\oa_8$:
\begin{alignat}{3}
\wt F_{3\tss 2}&=Y_{\al},\qquad
&&\wt F_{2\tss 1}=\wt F_{4\tss 3}=\wt F_{4'\tss 3}=\frac13\ts Y_{\be},
\qquad
&&\wt F_{3\tss 1}=-\wt F_{4\tss 2}=-\wt F_{4'\tss 2}=\frac13\ts Y_{\al+\be},
\non\\[0.4em]
\wt F_{4\tss 1}&=\wt F_{4'\tss 1}=\frac29\ts Y_{\al+2\tss\be},
\qquad &&\wt F_{3\pr\tss 1}=-\frac29\ts  Y_{\al+3\tss\be},\qquad
&&\wt F_{2\pr\tss 1}=\frac29\ts  Y_{2\tss\al+3\tss\be}.
\non
\end{alignat}
Furthermore,
\ben
\wt F_{1\tss 1}=-H_{\al}-\frac23\ts  H_{\be},\qquad
\wt F_{2\tss 2}=-H_{\al}-\frac13\ts H_{\be}, \qquad
\wt F_{3\tss 3}=-\frac13\ts H_{\be}\Fand \wt F_{4\tss 4}=0.
\een
Thus, the $9\times 9$ matrix introduced in the beginning of Sec.~\ref{sec:gend}
(with $n=4$) reduces to a matrix with entries in $\Vc(\p)$, where
$\p=\n_-\oplus\h$ is the subspace of $\g_2$.
The determinant of this matrix remains unchanged if the $(4,6)$ entry $-2\tss\di$
is replaced by $-\di$ and each entry of the lower left $4\times 4$ submatrix
is multiplied by $2$. Moreover, the determinant will also be unchanged
if we
replace row $6$ by the sum
of rows $4$ and $6$. Finally, the simultaneous expansion along the rows $5$ and $6$
reduces the determinant to the given form. This shows that
all elements $w_i$ belong to $\Wc(\g_2,f)$.

The centralizer $\g^{\ts f}_2$ is spanned by the elements $f$ and
$Y_{2\tss\al+3\tss\be}$.
The condition of Proposition~\ref{prop:genwalg} holds for the
elements $w_2,w_6$,
thus implying that they are generators
of the differential algebra $\Wc(\g_2,f)$.

The images of the elements $w_2,\dots,w_7$ under the homomorphism \eqref{phiw}
are the elements $\wt w_2,\dots,\wt w_7\in\Vc(\h)$ found from the relation
\begin{multline}
(\di+\wt F_{1\tss 1})(\di+\wt F_{2\tss 2})(\di+\wt F_{3\tss 3})\ts\di\ts
(\di-\wt F_{3\tss 3})(\di-\wt F_{2\tss 2})(\di-\wt F_{1\tss 1})\\[0.3em]
=\di^{\ts 7}+\wt w_2\ts \di^{\ts 5}+\wt w_3\ts \di^{\ts 4}+\wt w_4\ts \di^{\ts 3}
+\wt w_5\ts \di^{\ts 2}+\wt w_6\ts \di+\wt w_7.
\non
\end{multline}
In the notation of Sec.~\ref{subsec:ctt} we have
\ben
h^{(r)}_1=-H^{(r)}_{\al}\Fand
h^{(r)}_2=-H^{(r)}_{\be}.
\een
The differential algebra $\wt\Wc(\g_2,f)$ consists of the polynomials
in the variables $H^{(r)}_{\al}$ and $H^{(r)}_{\be}$, which are annihilated by
the screening operators $V_1$ and $V_2$ defined in \eqref{scr}.
Observe that the definition of the elements $\wt w_i$
coincides with that for the classical $\Wc$-algebra $\Wc(\oa_7,f)$
with the additional condition $\wt F_{2\tss 2}+\wt F_{3\tss 3}=\wt F_{1\tss 1}$;
see Sec.~\ref{sec:genb}. If this condition
is ignored, then the elements
$\wt w^{\tss(r)}_2,\wt w^{\tss(r)}_4$ and $\wt w^{\tss(r)}_6$
with $r$ running over nonnegative integers
are algebraically independent
generators of the algebra $\wt\Wc(\oa_{7},f)$; see \cite[Ch.~8]{f:lc}.
A direct calculation shows that the condition
$\wt F_{2\tss 2}+\wt F_{3\tss 3}=\wt F_{1\tss 1}$
implies
\ben
\wt w_4=\frac14\ts \wt w^{\ts 2}_2+3\tss \wt w^{\ts(2)}_2
\een
and that the family $\wt w^{\ts(r)}_2, \wt w^{\ts(r)}_6$ with $r\geqslant 0$ is algebraically
independent.
\epf

The injective homomorphism $\Wc(\g_{2},f)\hra \Vc(\h)$
taking $w_i$ to $\wt w_i$ can be regarded as
the {\it Miura transformation\/} associated with the Lie algebra $\g_2$;
see \cite[Ch.~8]{f:lc}.

\end{document}